\documentclass[12pt]{puthesis}
\usepackage{amsmath,amsthm,amsfonts,amssymb}

\theoremstyle{plain}
\newtheorem{theorem}{Theorem}[chapter]
\newtheorem{Cor}{Corollary}[chapter]
\newtheorem{Lemma}{Lemma}[chapter]
\newtheorem{proposition}{Proposition}[chapter]
\newtheorem*{GLH}{Generalized Lindel\"of Hypothesis}
\newtheorem*{CB}{Convexity Bound} 
\newtheorem*{SP}{Subconvexity Problem}

\theoremstyle{definition}
\newtheorem*{Def}{Definition}
\newtheorem*{Hyp}{Hypothesis}
\newtheorem{Remark}{Remark}[chapter]

\theoremstyle{remark}
\newtheorem{Case}{Case}
\newtheorem{Case2}{Case}
\newtheorem{Case3}{Case}
\newtheorem{Part}{Part}
\newtheorem*{Not}{Notation}

\DeclareMathOperator{\vol}{vol}
\DeclareMathOperator{\res}{res}

\newcommand{\csillag}{\sideset{}{^*}}
\newcommand{\smod}[1]{\text{\rm{ (mod $#1$)}}}
\newcommand{\nyil}{\rightsquigarrow}
\newcommand{\pet}{\langle\vfi,\vfi\rangle}

\newcommand{\jb}{\bigstar}
\newcommand{\ov}{\bar}

\newcommand{\GL}{{\rm GL}}
\newcommand{\SL}{{\rm SL}}

\newcommand{\lfi}{\lambda_\phi}
\newcommand{\lpsi}{\lambda_\psi}
\newcommand{\blfi}{{\bar\lambda}_\phi}
\newcommand{\rfi}{\rho_\phi}
\newcommand{\brfi}{{\bar\rho}_\phi}
\newcommand{\tmu}{{|\tilde\mu|}}

\newcommand{\al}{\alpha}
\newcommand{\be}{\beta}
\newcommand{\de}{\delta}
\newcommand{\dd}{\delta}
\newcommand{\ep}{\epsilon}
\newcommand{\gG}{\gamma}
\newcommand{\GG}{\Gamma}
\newcommand{\ka}{\kappa}
\newcommand{\la}{\lambda}
\newcommand{\lL}{\lambda}
\newcommand{\LL}{\Lambda}
\newcommand{\om}{\omega}
\newcommand{\si}{\sigma}
\newcommand{\ff}{\phi}
\newcommand{\vfi}{\phi}

\newcommand{\cc}{\mathfrak{c}}
\newcommand{\FI}{\mathfrak{I}}
\newcommand{\II}{\mathfrak{I}}

\newcommand{\cC}{\mathcal{C}}
\newcommand{\FF}{\mathcal{F}}
\newcommand{\HH}{\mathcal{H}}
\newcommand{\CQ}{\mathcal{Q}}

\newcommand{\QQ}{\mathbb{Q}}
\newcommand{\RR}{\mathbb{R}}
\newcommand{\CC}{\mathbb{C}}
\newcommand{\ZZ}{\mathbb{Z}}

\author{Gergely Harcos}
\title{New Bounds for Automorphic $L$-functions}
\submitted{June, 2003}

\abstract{This dissertation contributes to the analytic theory of automorphic
$L$-functions.

We prove an approximate functional equation for the central value
of the $L$-series attached to
an irreducible cuspidal automorphic representation $\pi$ of $\GL_m$ over
a number field.
The approximation involves a smooth truncation of the Dirichlet series
$L(s,\pi)$ and $L(s,\tilde\pi)$ after about $\sqrt{C}$ terms,
where $C$ denotes the analytic conductor 
(of $\pi$ and $\tilde\pi$ at the central point) 
introduced by Iwaniec and Sarnak.
We investigate the decay rate of the cutoff function and its derivatives. 
We also see
that the truncation can be made uniformly explicit at the
cost of an error term.
The results extend to products of central values.

We establish, via the Hardy--Littlewood circle method,
a nontrivial bound on shifted convolution sums
of Fourier coefficients coming from classical holomorphic or Maass
cusp forms of arbitrary level and nebentypus. These sums are
analogous to the binary additive divisor sum which has been
studied extensively. We achieve polynomial uniformity in all the 
parameters of the cusp forms by carefully estimating the Bessel
functions that enter the analysis.
As an application we derive, extending work
of Duke, Friedlander and Iwaniec, a subconvex estimate on the
critical line for $L$-functions associated to character twists of
these cusp forms.

We also study the shifted convolution sums via the Sarnak--Selberg spectral
method. For holomorphic cusp forms this approach detects optimal cancellation 
over any totally real number field. For Maass cusp forms the method is
burdened with complicated integral transforms. We succeed in inverting the simplest
of these transforms whose kernel is built up of Gauss hypergeometric functions.

\medskip
\noindent}
\acknowledgements{My gratitude to my advisor, Peter Sarnak, is due not only to his guidance
in the preparation of the current work, but, most of all, to his permanent
influence on my mathematical outlook. I feel fortunate to have glimpsed,
at times, the mathematical purview provided by his wide perspective.

I should also manifest my debt to all my teachers, and to all companions
who helped to form me as a mathematician. My friends, in general, are to
be thanked for their valued support. What I owe to my family cannot
be distinguished from who I am and who I am becoming; I can only ascertain
that most of my worth comes from them. Before all, I thank my father, 
P\'eter Harcos, and Yvette Vajda, my bride.

Copyright notice. Chapter~\ref{functional} first appeared, in almost identical form,
as ``Uniform approximate functional equation for principal $L$-functions,''
\emph{Int. Math. Res. Not.}
\textbf{2002}, no. 18,
923--932, published by Hindawi Publishing Corporation. This part of the work may be reproduced by any means for educational and scientific purposes without fee or permission with the exception of reproduction by services that collect fees for delivery of documents.
}
\dedication{To My Father}

\begin{document}
\chapter{Introduction}

\section{Prologue}

$L$-functions are among the most fundamental and most fascinating objects in number theory.
An $L$-function can be attached to
\begin{itemize}
\item[(1)] a smooth projective variety defined over a number field (Hasse, Weil),
\item[(2)] an irreducible complex or $l$-adic representation of the Galois group 
of a number field (Artin, Grothendieck), or
\item[(3)] a cusp form or irreducible cuspidal automorphic representation (Hecke, Langlands, Godement--Jacquet).
\end{itemize}
An $L$-function is defined in terms of local data. In each of the cases
above, this local data consists of
\begin{itemize}
\item[(1)] the number of points of the reduction of the projective variety
to various finite fields,
\item[(2)] the eigenvalues of the Frobenius elements in the Galois group, or
\item[(3)] the Langlands parameters of the automorphic form or representation.
\end{itemize}
By definition, the $L$-function is given as an Euler product over the
rational primes of the local data:
\[L(s)=\prod_p L_p(s).\]
Various results and conjectures relating these objects
add up to the general philosophy that every $L$-function of arithmetic nature is
a ratio of automorphic $L$-functions.

Besides their many combinatorial and algebraic properties, $L$-functions are very much
analytic objects. 
Understanding their analytic behaviour is an important task, especially if it gives rise
to arithmetic implications. A classical example is Chebotarev's density theorem on Frobenius
elements in the Galois group. The analytic properties of an $L$-function
are most accessible when the $L$-function
is known to come from an automorphic form. Even in this case, our knowledge is surprisingly limited.
It has been realized only recently how widely such knowledge could be applied to deep diophantine problems.

\section{Size of an $L$-function}\label{size}

A foremost issue in such applications is that of the size of an $L$-function.
Let us first fix our notation for a general discussion. 
We consider a number field $F$
of degree $d$ and an irreducible cuspidal automorphic representation 
$\pi$ of $\GL_m$ over $F$ with unitary central character. 
By Flath's theorem, $\pi$ can be written uniquely as a restricted tensor product $\otimes_v\pi_v$, where $\pi_v$ is an irreducible admissible representation  of $\GL_m(F_v)$ for each place $v$ of $F$.
Accordingly, the complete $L$-function associated to $\pi$ is defined as a product of
local $L$-functions,
\[\LL(s,\pi)=\prod_v L(s,\pi_v).\]
It is convenient to collect the local factors for $v$ underlying a given rational place $w$,
and introduce the subproducts
\[L(s,\pi_w)=\prod_{v\mid w}L(s,\pi_v).\]
For the infinite place $w=\infty$ the subproduct takes the form
\begin{equation}\label{eq81}
L(s,\pi_\infty)=\prod_{j=1}^{md}\pi^{\frac{\mu_j-s}{2}}\GG\left(\frac{s-\mu_j}{2}\right),
\end{equation}
while for a finite rational prime $w=p$ we have
\begin{equation}\label{eq72}
L(s,\pi_p)=\prod_{j=1}^{md}\frac{1}{1-\al_j(p)p^{-s}}.
\end{equation}
(Note that $\pi$ 
inside the first product refers to the positive constant, not the representation.)
The numbers $\mu_j$ (resp. $\al_j(p)$) are called the Archimedean (resp. non-Archimedean)
Langlands parameters and satify the following uniform bound by Theorem~1 of \cite{LRS}.

\begin{theorem}[Luo--Rudnick--Sarnak]\label{ThLRS}
\begin{equation}\label{eq70}
\sup\{\Re\mu_j,\Re\log_p\al_j(p)\}\leq\frac{1}{2}-\frac{1}{m^2+1}.
\end{equation}
\end{theorem}

The local factor $L(s,\pi_\infty)$ is distinguished in the sense that in vertical strips
it decays exponentially while the other factors $L(s,\pi_p)$ remain bounded away from $0$. 
This fact alone provides
ample justification for isolating the finite part
\begin{equation}\label{eq73}
L(s,\pi)=\prod_{p<\infty}L(s,\pi_p),\quad\Re s>\frac{3}{2}-\frac{1}{m^2+1},
\end{equation}
an absolutely convergent Euler product over the rational primes
by (\ref{eq70}). The resulting complete
\[\LL(s,\pi)=L(s,\pi_\infty)L(s,\pi)\]
extends to an entire function which is bounded in vertical strips 
(except for $\pi=|\det|^{it}$ when a 
simple pole occurs at $s=1-it$), and satisfies a functional equation of the form
\begin{equation}\label{eq71}
N^\frac{s}{2}\LL(s,\pi)=
\ka N^\frac{1-s}{2}\LL(1-s,\tilde\pi).
\end{equation}
$N$ is the arithmetic conductor (a positive integer), $\ka$ is the root
number (of modulus 1), and $\tilde\pi$ is the contragradient representation of $\pi$.
The local $L$-functions of $\pi$ and $\tilde\pi$ are connected by
\begin{equation}\label{eq82}
\ov L(s,\pi_v)=L(\ov{s},\tilde\pi_v).
\end{equation}

It is natural to expect that 
$L(s,\pi)$ has a moderate size in vertical strips, so that
$\LL(s,\pi)$ inherits the exponential decay of the Archimedean factor $L(s,\pi_\infty)$.
We shall formulate a more precise and more general statement
using the analytic conductor introduced by Iwaniec and Sarnak \cite{IS}:
\[C(s,\pi)=\frac{N}{(2\pi)^{md}}\prod_{j=1}^{md}|s-\mu_j|.\]
In order to bound automorphic $L$-functions, it is essential to 
represent them as absolutely convergent Dirichlet series
\begin{equation}\label{eq76}
L(s,\pi)=\sum_{n=1}^\infty\frac{\la_\pi(n)}{n^s}.
\end{equation}
Certainly (\ref{eq72}), (\ref{eq70}) and (\ref{eq73}) guarantee 
that $L(s,\pi)$ acquires this
form in the half-plane $\Re s>\frac{3}{2}-\frac{1}{m^2+1}$. As a
by-product, we also see that the coefficients satisfy
\[\la_\pi(n)\ll_{\ep,m,d}n^{\frac{1}{2}-\frac{1}{m^2+1}+\ep}\]
for any $\ep>0$. Upon the Ramanujan--Selberg conjectures we could
replace the occurrences of $\frac{1}{2}-\frac{1}{m^2+1}$ in (\ref{eq70})
and in the previous inequality by 0. These improved bounds hold unconditionally
in a certain average form by Theorem 4 of \cite{Mol}.

\begin{theorem}[Molteni]\label{ThMol} Uniformly in $\ep>0$ and $x>0$,
\begin{equation}\label{eq74}
\sum_{n\leq x}|\la_\pi(n)|\ll_\ep x^{1+\ep}C\left(\tfrac{1}{2},\pi\right)^\ep.
\end{equation}
The implied constant depends only on $\ep$, $m$ and $d$.
\end{theorem}

It should be noted that Molteni assumes $\Re\mu_j\leq 0$ for all $j$
(cf. axiom (A4) in \cite{Mol}), but his argument works equally
well with the weaker bound (\ref{eq70}). In particular, the Dirichlet series 
(\ref{eq76}) is absolutely
convergent in the larger half-plane $\Re s>1$, and it satisfies
\begin{equation}\label{eq77}
L(\si,\pi)\ll_{\si,\ep,m,d}C\left(\tfrac{1}{2},\pi\right)^\ep,\quad\si>1.
\end{equation}
By replacing $\pi$ with its twist $\pi\otimes|\det|^{it}$ this becomes
\begin{equation}\label{eq78}
L(\si+it,\pi)\ll_{\si,\ep,m,d}C\left(\tfrac{1}{2}+it,\pi\right)^\ep,\quad\si>1.
\end{equation}
We can combine (\ref{eq77}) with the functional equation (\ref{eq71})
to deduce uniform bounds in the half-plane $\Re s<0$. First,
\[L(\si,\pi)\ll_{\si,\ep,m,d}C\left(\tfrac{1}{2},\pi\right)^{1/2-\si+\ep},\quad\si<0,\]
and then, by replacing $\pi$ with $\pi\otimes|\det|^{it}$,
\begin{equation}\label{eq79}
L(\si+it,\pi)\ll_{\si,\ep,m,d}C\left(\tfrac{1}{2}+it,\pi\right)^{1/2-\si+\ep},\quad\si<0.
\end{equation}
Finally, we can interpolate between (\ref{eq78}) and (\ref{eq79}) by
the Phragm\'en--Lindel\"of convexity principle to obtain bounds inside
the critical strip $0<\Re s<1$ (or on the boundaries away from the possible
pole).

\begin{CB}  For any $0<\si<1$ and any $\ep>0$, there is a uniform bound
\begin{equation}\label{eq80}
L(\si+it,\pi)\ll_{\si,\ep}C\left(\tfrac{1}{2}+it,\pi\right)^{(1-\si)/2+\ep}.
\end{equation}
The implied constant depends only on $\si$, $\ep$, $m$ and $d$.
\end{CB}

The expontents given by (\ref{eq78}) and (\ref{eq79}) are sharp.
We expect, however, that a much stronger inequality holds in place of the
convexity bound. 

\begin{GLH} For any $0<\si<1$ and any $\ep>0$, there is a uniform bound
\begin{equation}\label{eq84}
L(\si+it,\pi)\ll_{\si,\ep}
C\left(\tfrac{1}{2}+it,\pi\right)^{\max(0,1-2\si)/2+\ep}.
\end{equation}
The implied constant depends only on $\si$, $\ep$, $m$ and $d$.
\end{GLH}

This very powerful statement is a consequence of the
generalized Riemann hypothesis that all the roots of $\LL(s,\pi)$ lie
on the critical line $\Re s=\frac{1}{2}$.
In fact, the resolution of several deep equidistribution questions in number theory
relies on a small but substantial improvement on the convexity bound in
certain families of automorphic $L$-functions. For convenience and applicability
we focus on the critical line $\Re s=\frac{1}{2}$.

\begin{SP} Show that there is a $\dd=\dd(m,d)>0$ such that
\begin{equation}\label{eq83}
L(s,\pi)\ll_{m,d} C\left(s,\pi\right)^{1/4-\dd}, \quad\Re s=\frac{1}{2}.
\end{equation}
\end{SP}

Applications include equidistribution of lattice points on ellipsoids 
(Linnik's problem), characterization of integers represented by a given 
quadratic form over a number field
(Hilbert's 11th problem), equidistribution of certain 
Galois orbits of CM-points on
Shimura varieties (evidence toward the Andr\'e--Oort conjecture), 
and equidistribution of mass in arithmetic quantum chaos.

\section{Approximate functional equation}\label{functional_intro}

It is not obvious that the coefficients $\la_\pi(n)$ can be used to reveal
the finer behaviour of $L(s,\pi)$ in the critical strip $0<\Re s<1$. This was 
originally realized for the Riemann zeta function 
\[\zeta(s)=\sum_{n=1}^\infty\frac{1}{n^s}\]
by Hardy and Littlewood in 1921 \cite{HL2}. They established an approximation to $\zeta(s)$,
called an \emph{approximate functional equation}, a special case of which
reads as follows:
\[\zeta\left(\frac{1}{2}+it\right)=\sum_{n\leq\sqrt{\frac{|t|}{2\pi}}}\frac{1}{n^{\frac{1}{2}+it}}+
\frac{\zeta\left(\frac{1}{2}+it\right)}{\zeta\left(\frac{1}{2}-it\right)}
\sum_{n\leq\sqrt{\frac{|t|}{2\pi}}}\frac{1}{n^{\frac{1}{2}-it}}+O\bigl(|t|^{-\frac{1}{4}}\log|t|\bigr).\]
Note that the factor in front of the second sum is of modulus 1 and does not destroy
the symmetry $t\leftrightarrow -t$. This formula was extended and studied by many researchers with
focus generally restricted to small powers of Dirichlet $L$-functions or Dedekind $L$-functions.
Among the few studies with a larger scope the most notable ones are by
Chandrasekharan and Narasimhan \cite{CN}, Lavrik \cite{La}, and Ivi\'c \cite{Iv}.

In Chapter~\ref{functional} we shall present uniform variants of the approximate functional equation for all automorphic $L$-functions. We shall demonstrate that the values of
$L(s,\pi)$ on the critical line $\Re s=\frac{1}{2}$ can be approximated as a sum of
two Dirichlet series which have essentially $\sqrt{C(s,\pi)}$ terms. The relevance of the analytic conductor has not been displayed in this general context before. 
In fact, we had to do some ``fine tuning'' on the original analytic conductor
of Iwaniec and Sarnak \cite{IS} in order to achieve our goal.

The result we obtain fits well into the philosophy that $L$-functions (or rather, $L$-values) should be considered in families \cite{IS}. We shall employ smooth cutoff functions as they are more natural for the problem and also yield better error terms.
First we obtain an exact representation by an implicit cutoff function
with uniform decay properties (Theorem~\ref{theorem01}). This
formula is most useful for families whose Archimedean parameters
remain bounded. The second representation (Theorem~\ref{theorem02}),
inspired by the recent work of Ivi\'c \cite{Iv}, has a more explicit 
main term at the cost of an error term. This formula
works best in families where the Archimedean parameters grow large simultaneously. The proofs are based on standard Mellin transform techniques, and they make crucial use of the estimates of Luo--Rudnick--Sarnak (\ref{eq70}) and
Molteni (\ref{eq74}). A variant of the method yields
similar formulae for products of central values (e.g. for higher
moments).

\section{Amplification}\label{amplification}

The approximate functional equation reduces the subconvexity problem to cancellation in finite smooth sums
\[S(X,\pi)=\sum_{n=1}^\infty\la_\pi(n)w\left(\frac{n}{X}\right),\]
where $w:(0,\infty)\to\CC$ is a fixed weight function of compact support on the positive axis. More precisely, by combining Corollary~\ref{Cor1} with a smooth decomposition of unity, we can see that a variant of (\ref{eq83}), 
\begin{equation}\label{eq83b}
\forall\ep>0:\forall t\in\RR:\quad L(\tfrac{1}{2}+it,\pi)\ll_{\ep,m,d} C\left(\tfrac{1}{2}+it,\pi\right)^{1/4-\dd+\ep}, \end{equation}
follows from a uniform bound\begin{equation}\label{eq85}
S(X,\pi)\ll_{w,m,d}C\left(\tfrac{1}{2},\pi\right)^{1/4-\dd+\ep}\sqrt{X}
\end{equation}
in the range $X\leq C\left(\tfrac{1}{2},\pi\right)^{1/2+\ep}$. It should be observed that Molteni's bound (\ref{eq74}) yields an even stronger estimate whenever $X\leq C\left(\tfrac{1}{2},\pi\right)^{1/2-2\dd}$. The above inequality (with no restriction on $X$) is
in fact equivalent to the subconvex bound (\ref{eq83b}), as can be seen from the representation
\[S(X,\pi)=\int_{1/2-i\infty}^{1/2+i\infty}L(s,\pi)X^s W(s)\,ds,\]where
\[W(s)=\int_0^\infty w(x)x^s\frac{dx}{x}\]
denotes the Mellin transform of $w(x)$. 

By this line of thought we also see that the generalized Lindel\"of hypothesis 
(\ref{eq84}) translates into strong square-root cancellation among the coefficients $\la_\pi(n)$:
\[S(X,\pi)\ll_{\ep,w,m,d}C\left(\tfrac{1}{2},\pi\right)^{\ep}\sqrt{X}.\]
In particular, we expect that in a family $\FF$ of cusp forms $\pi$ we have
\[\frac{1}{|\FF|}\sum_{\pi\in\FF}|S(X,\pi)|^2\ll_{\ep,w,m,d}C^\ep X,\]
as long as the analytic conductors satisfy $C\left(\tfrac{1}{2},\pi\right)\asymp C$. It is often possible to apply ideas from harmonic analysis to establish the preceding square mean bound for certain families $\FF$. As an immediate consequence, we obtain a pointwise bound
\[S(X,\pi)\ll_{\ep,w,m,d}C^\ep\sqrt{|\FF|X},\quad\pi\in\FF.\]
If we can guarantee that $|\FF|\ll C^{1/2-2\dd}$, then a subconvex bound for $L\left(\frac{1}{2},\pi\right)$ is established in the form (\ref{eq85}).
In most cases, however, harmonic analysis just falls short of establishing subconvexity. This is not surprising in the light of the extensive deep applications of subconvex bounds in number theory. The roots of subconvexity lie in arithmetic.

Amplification is an arithmetic device to substitute for shortening the family $\FF$. It appeared in the seminal work of Duke, Friedlander and Iwaniec
\cite{FI,DFI2}. The basic idea is to introduce nonnegative arithmetic weights $|a_\pi|^2$ so that
\[\frac{1}{|\FF|}\sum_{\pi\in\FF}|a_\pi|^2|S(X,\pi)|^2\ll_{\ep,w,m,d}C^\ep X,\]
while $|a_\pi|$ is larger than $C^\dd$ for a specific $\pi\in\FF$ and some $\dd>0$. Then we only need to guarantee that $|\FF|\ll C^{1/2+\ep}$, and subconvexity follows. The details in carrying out this program can become very complicated. Much of this thesis is devoted to study shifted convolution sums of the coefficients $\la_\pi(n)$, the sums that lie at the heart of the amplification method in the cases where it is known to work.

\section{Shifted convolution sums and the circle method}\label{shifted}

A particularly interesting (conjectural) family of automorphic representations consists of Rankin--Selberg products $\pi\otimes\rho$, where $\pi$ is a fixed cusp form on $\GL_m$ and $\rho$ varies over cusp forms on a fixed $\GL_n$ ($n\leq m$). The $L$-functions $L(s,\pi\otimes\rho)$ can be defined intrinsically and the expected analytic properties have been established by the work of many authors. The approach of amplification to establish subconvexity for these $L$-functions naturally leads to shifted convolution sums for $\pi$:
\begin{equation}\label{eq100}
D_f(a,b;h)=\sum_{am\pm bn=h}\la_\pi(m)\ov\la_\pi(n)f(am,bn).
\end{equation}
Here $a$, $b$, $h$ are positive integers
and $f$ is some nice weight function on $(0,\infty)\times (0,\infty)$,
e.g. smooth and
compactly supported on a box $[X,2X]\times[Y,2Y]$. If we have a uniform estimate
\[\sum_{m\leq x}|\la_\pi(m)|^2\ll_\pi x,\]
then the size of the sum (\ref{eq100}) can be seen to be
at most $O_{\ep,f,\pi}\bigl(\sqrt{XY}\bigr)$.
In order to achieve subconvexity, we  need
to improve on this bound in the $X$ and $Y$
aspects with certain uniformity regarding the other parameters.

Historically, the first examples of shifted convolution sums were
generalized binary additive divisor sums, whose coefficients are given
in terms of the divisor function:
\[D_f^\tau(a,b;h)=\sum_{am\pm bn=h} \tau(m)\tau(n)f(am,bn).\]
Note that the $\tau(n)$'s generate
$\zeta^2(s)$, and they also appear as Fourier coefficients of the modular form
$\frac{\partial}{\partial s}E(z,s)\big|_{s=1/2}$, where $E(z,s)$ is
the Eisenstein series for $SL_2(\ZZ)$. 
These sums have been studied extensiviely since 1926, when Kloosterman published his famous refinement of the circle method \cite{Kl}. A short summary of subsequent developments can be found in \cite{DFI}.

The crucial insight of Kloosterman was to make use of the very regular distribution of Farey fractions on the unit interval. By applying Vorono\"\i-type summation formulae for the relevant exponential generating functions (which in turn reflect modular transformation properties), the binary additive sum in question decomposes to a main term and an error term in a natural fashion. The main term arises, because $E(z,s)$ is not cuspidal, and the error term is expressed in terms of Kloosterman sums
\[S(m,n;q)=\csillag\sum_{d\smod{q}}e_q\bigl(dm+{\bar
d}n\bigr),\]
for which a nontrivial bound is needed. Kloosterman \cite{Kl} did provide a nontrivial bound, and later Weil \cite{W} and Esterman \cite{E} proved the optimal estimate. 

This classical approach was revived recently by Duke, Friedlander and Iwaniec \cite{DFI}. The Farey dissection being disguised as the $\dd$-method, the Vorono\"\i-type summation formula is still utilized at all frequencies so as to yield the following general result.

\begin{theorem}[Duke--Friedlander--Iwaniec]
Let $a$, $b$ coprime and assume that the partial derivatives of the
weight function $f$ satisfy the estimate
\begin{equation}\label{eq1}
x^ky^lf^{(k,l)}(x,y)\ll_{k,l}
\left(1+\frac{x}{X}\right)^{-1}\left(1+\frac{y}{Y}\right)^{-1}P^{k+l}
\end{equation}
with some $P,X,Y\geq 1$ for all $k,l\geq 0$. Then
\[D_f^\tau(a,b;h)=\int_0^\infty g(x,\mp x\pm h)\,dx+
O\bigl(P^{5/4}(X+Y)^{1/4}(XY)^{1/4+\ep}\bigr),\]
where the implied constant depends only on $\ep$,
\[g(x,y)=f(x,y)\sum_{q=1}^\infty \frac{(ab,q)}{abq^2}c_q(h)
(\log x-\lL_{aq})(\log y-\lL_{bq}),\]
$c_q(h)=S(h,0;q)$ denotes Ramanujan's sum, and
$\lL_{aq}$, $\lL_{bq}$ are constants given by
\[\lL_{aq}=2\gG+\log\frac{aq^2}{(a,q)^2}.\]
\end{theorem}

As was pointed out in \cite{DFI}, the error term is smaller than
the main term whenever
\[P^{5/4}ab\ll(X+Y)^{-5/4}(XY)^{3/4-\ep}.\]

In Chapter~\ref{maass} we shall extend the above ideas
to exhibit nontrivial cancellation in the shifted convolution sums (\ref{eq100}) for cuspidal automorphic representations $\pi$ of $\GL_2$ over $\QQ$.
 In fact, we shall estimate the more general sums
\begin{equation}\label{eq24}
D_f(a,b;h)=\sum_{am\pm bn=h} \lfi(m)\lpsi(n)f(am,bn),
\end{equation}
where $\lfi(m)$ (resp. $\lpsi(n)$) are the normalized Fourier coefficients of a
classical holomorphic or weight zero Maass cusp form $\phi$ (resp. $\psi$) of 
arbitrary level and nebentypus. The conclusion is recorded in Theorem~\ref{Th1}. 
In Chapter~\ref{subconvex} we shall apply the result about shifted convolution sums to obtain a subconvex bound for the values
$L(s,\vfi\otimes\chi)$, where $\vfi$ is a primitive form in the sense of Atkin--Lehner theory \cite{AL,Li,ALi}, $s$ is a fixed point on the critical line, and
$\chi$ runs through primitive Dirichlet characters of conductor prime to the level of $\vfi$ (Theorem~\ref{Th2}). A specialization to the central point $s=1/2$ yields, via Waldspurger's theorem and its
generalization \cite{Wa,Sh}, nontrivial bounds for the Fourier-coefficients
of holomorphic or Maass cusp forms of half-integral weight. These bounds in turn can be applied to resolve 
Linnik's problem \cite{D,D-SP}.

\section{Shifted convolution sums and spectral theory}\label{shifted2}

A different spectral approach was developed by Sarnak for all levels. The method can be traced back to the discovery of Rankin and Selberg, that for a holomorphic cusp form
\[\phi(z)=\sum_{n=1}^\infty\rfi(n)e(nz)\]
of weight $k$, level $N$ and arbitrary nebentypus, there is an integral representation
\begin{equation}\label{eq300}
\sum_{n=1}^\infty\frac{|\rfi(n)|^2}{n^{s+k-1}}=\frac{(4\pi)^{s+k-1}}{\GG(s+k-1)}
\int\limits_{\GG\backslash\HH}y^k|\phi(z)|^2E(z,s)\frac{\,dx\,dy}{y^2},
\end{equation}
where $\GG\backslash\HH$ is a fundamental domain for the action of the Hecke congruence subgroup $\GG=\GG_0(N)$ on the upper half-plane $\HH=\{x+iy:y>0\}$, and 
\[E(z,s)=\sum_{\gG\in\GG_\infty\backslash\GG}y^s(\gG z)\]
denotes Eisenstein's series. The above identity can be proved by a simple unfolding technique, and it shows that the summatory function of the coefficients $|\rfi(n)|^2$ depends largely on the analytic properties of $E(z,s)$. The Eisenstein series $E(z,s)$ 
is a meromorphic function in the $s$-plane with only finitely many poles 
in $\Re s\geq 1/2$. There is always a simple pole at $s=1$ with residue explicitly given by
\[\res\displaylimits_{s=1}E(z,s)=\frac{1}{\vol\bigl(\GG\backslash\HH\bigr)}.\]
The other poles are also simple and lie on the real segment $(1/2,1)$.
They correspond to the residual spectrum of $\GG\backslash\HH$, a set which was anticipated by Selberg to be empty.

The connection with the shifted convolution sums (\ref{eq100}) becomes apparent if we specify $\GG=\GG_0(Nab)$, 
replace $E(z,s)$ by the Poincar\'e series
\[P_h(z,s)=\sum_{\gG\in\GG_\infty\backslash\GG}y^s(\gG z)e(-hx(\gG
z)),\]
and the $\GG_0(N)$-invariant product $y^k|\phi(z)|^2$ by the $\GG$-invariant product $y^k\phi(az)\ov\phi(bz)$. We obtain, by the same unfolding technique,
\begin{equation}\label{eq310}
\sum_{am-bn=h}\frac{\rfi(m)\brfi(n)}{(am+bn)^{s+k-1}}=\frac{(2\pi)^{s+k-1}}{\GG(s+k-1)}
\int\limits_{\GG\backslash\HH}y^k\phi(az)\ov\phi(bz)P_h(z,s)\frac{\,dx\,dy}{y^2}.\end{equation}
The integral equals, by definition, the Petersson inner product of the $\GG$-invariant functions $U(z)=y^k\phi(az)\ov\phi(bz)$ and $\ov P_h(z,s)$, and it
can be decomposed according
to the spectrum of $L^2(\GG\backslash\HH)$. The discrete part of
the spectrum corresponds to an orthonormal basis of Maass cusp forms
\[\begin{split}
\ff_0(x+iy)&=\frac{1}{\vol^{1/2}(\GG\backslash\HH)},\\\\
\ff_j(x+iy)&=\sqrt{y}\sum_{n\neq 0}\la_j(n)K_{i\tau_j}\bigl(2\pi|n|y\bigr)e(nx),
\quad j=1,2,\dots,
\end{split}\]
while the continuous spectrum is provided by
the Eisenstein series 
\[E_\cc(\cdot,\tfrac{1}{2}+i\tau)=\delta_\cc y^s+\eta_\cc(s)y^{1-s}+
\sqrt{y}\sum_{n\neq 0}\la_{\cc,\tau}(n)K_{i\tau}\bigl(2\pi|n|y\bigr)e(nx)
,\quad\tau\in\RR,\]
where $\cc$ is a singular cusp of $\GG\backslash\HH$.
The decomposition 
reads, at least formally, as
\begin{multline*}I(s)=\langle U,\ov P_h(.,s)\rangle=\sum_{j=0}^\infty
\langle U,\ff_j\rangle\langle \ff_j,\ov P_h(\cdot,s)\rangle\\
+\sum_\cc\frac{1}{4\pi}\int_{-\infty}^\infty 
\langle U,E_\cc(\cdot,\tfrac{1}{2}+i\tau)\rangle
\langle E_\cc(\cdot,\tfrac{1}{2}+i\tau),\ov P_h(\cdot,s)\rangle
\,d\tau.\end{multline*} 
We observe that the inner products $\langle\ff_j,\ov P_h(.,s)\rangle$
and $\langle E_\cc(\cdot,\tfrac{1}{2}+i\tau),\ov P_h(\cdot,s)\rangle$ can be unfolded to
\[\begin{split}
\langle\ff_j,\ov P_h(.,s)\rangle
&=\frac{\la_j(h)}{4(\pi h)^{s-\frac{1}{2}}}
\GG\left(\frac{s-\frac{1}{2}+i\tau_j}{2}\right)
\GG\left(\frac{s-\frac{1}{2}-i\tau_j}{2}\right),\\\\
\langle E_\cc(\cdot,\tfrac{1}{2}+i\tau),\ov P_h(\cdot,s)\rangle
&=\frac{\la_{\cc,\tau}(h)}{4(\pi h)^{s-\frac{1}{2}}}
\GG\left(\frac{s-\frac{1}{2}+i\tau}{2}\right)
\GG\left(\frac{s-\frac{1}{2}-i\tau}{2}\right),
\end{split}\]
where $\tfrac{1}{4}+\tau_j^2$ (resp. $\tfrac{1}{4}+\tau^2$) 
denotes the Laplacian eigenvalue of $\ff_j$ (resp. of
$E_\cc(\cdot,\tfrac{1}{2}+i\tau)$).

It follows that the size of $I(s)$ (including the location of its poles) are determined by the residual spectrum of $\GG\backslash\HH$, the size of the Fourier coefficients $\la_j(h)$ and $\la_{\cc,\tau}(h)$, and the size of the triple products $\langle U,\ff_j\rangle$ and $\langle U,E_\cc(\cdot,\tfrac{1}{2}+i\tau)\rangle$. We know that $\la_{\cc,\tau}(h)$ is of size at most $h^\ep$, and the Ramanujan conjecture predicts the same for $\la_j(h)$. In addition, the Selberg conjecture predicts that the residual spectrum is empty, that is, $\tau_j\in\RR$.
As a substitute for these conjectures, we shall only assume the following statement which is known for many nontrivial values $\theta<1/2$ (cf. (\ref{eq70})):

\begin{Hyp}For any cusp form $\pi$ on $\GL_2$ over $\QQ$, the local Langlands
parameters $\mu_{j,\pi}$ and $\al_{j,\pi}(p)$ $(j=1,2)$ satisfy
\[\begin{split}
|\Re\mu_{j,\pi}|\leq\theta,&\quad\text{if $\pi_\infty$ is unramified;}\\
\bigl|\Re\log_p\al_{j,\pi}(p)\bigr|\leq\theta,&\quad\text{if $\pi_p$ is unramified\  ($p<\infty$).}
\end{split}\]
\end{Hyp}
The behaviour of the triple products $\langle U,\ff_j\rangle$ and $\langle U,E_\cc(\cdot,\tfrac{1}{2}+i\tau)\rangle$ was only understood recently by Sarnak \cite{Sa4,Sa2}. He showed that
\[\langle U,\ff_j\rangle\ll_{\phi}\bigl(1+|\tau_j|\bigr)^{k+1}e^{-\frac{\pi}{2}|\tau_j|},\]
and similarly for $\langle U,E_\cc(\cdot,\tfrac{1}{2}+i\tau)\rangle$. Note that the exponential decay in the eigenvalue parameter $\tau_j$ (resp. $\tau$) exactly compensates the exponential decay of the coefficient $\GG(s+k-1)$ in (\ref{eq310}).

If $\ff_1,\ff_2,\dots$ are suitably chosen Maass--Hecke
cuspidal eigenforms, then this argument leads to the powerful estimate
\begin{equation}\label{eq311}
J(s)=\sum_{am-bn=h}\frac{\rfi(m)\brfi(n)}{(am+bn)^{s+k-1}}\ll_{\phi,\ep}
(ab)^{1-\frac{k}{2}}h^{\frac{1}{2}+\theta-\si+\ep}|s|^3,\quad\Re s\geq\frac{1}{2}+\theta+\ep.
\end{equation} 
Note that $\theta=7/64$ is eligible by the recent work of Kim and Sarnak \cite{Ki}. The strength of this result comes from the fact that it can be combined with the technique of Mellin transforms to yield a nontrivial bound for any shifted convolution sum
\[\sum_{am-bn=h}\lfi(m)\blfi(n)W\left(\frac{am+bn}{h}\right),\]
where $W$ is an arbitrary smooth function $(1,\infty)\to\CC$ of compact support, and
\[\lfi(m)=m^{\frac{1-k}{2}}\rfi(m)\]
denotes the normalized Fourier coefficients of $\phi$. To see this connection, we introduce for convenience the variable
\[u=\frac{am+bn}{h},\]
as well the function
\[V(u)=(1-u^{-2})^{\frac{1-k}{2}}W(u),\]
then for any $\si>1$ we get
\[\begin{split}\sum_{am-bn=h}\lfi(m)\blfi(n)W(u)
&=(4ab)^{\frac{k-1}{2}}\sum_{am-bn=h}\frac{\rfi(m)\brfi(n)}{(am+bn)^{k-1}}V(u)\\\\
&=\frac{1}{2\pi i}\int_{(\si)}(4ab)^{\frac{k-1}{2}}h^sJ(s)\hat V(s)\,ds.
\end{split}\]
We can rewrite (\ref{eq311}) as
\[(4ab)^{\frac{k-1}{2}}h^sJ(s)\ll_{\phi,\ep}
(ab)^{\frac{1}{2}}h^{\frac{1}{2}+\theta+\ep}|s|^3,\quad\Re s\geq\frac{1}{2}+\theta+\ep,\]
therefore by shifting $\si>1$ to any $\si>\frac{1}{2}+\theta$ we can conclude that
\[\sum_{am-bn=h}\lfi(m)\blfi(n)W(u)\ll_{\phi,\ep}
(ab)^{\frac{1}{2}}h^{\frac{1}{2}+\theta+\ep}\sup_{\si+i\RR}\bigl|s^3\hat V(s)\bigr|.\]In particular, if $W$ is supported on $(X,2X)$, then we obtain
\[\sum_{am-bn=h}\lfi(m)\blfi(n)W(u)\ll_{\phi,\si,\ep}
(ab)^{\frac{1}{2}}h^{\frac{1}{2}+\theta+\ep}X^{\si}\max_{j=0,1,2,3}\bigl\|V^{(j)}\bigr\|_\infty,
\quad\si\geq\frac{1}{2}+\theta+\ep.\]

For a Maass cusp form of weight $\ka$ and level $N$ the analogous argument leads to complicated integral transforms. For such a form $\phi$ the Fourier expansion reads
\[\phi(x+iy)=\sum_{n\neq 0}
\rfi(n)\tilde W_{\frac{n}{|n|}\frac{\ka}{2},i\mu}\bigl(4\pi|n|y\bigr)e(nx),\]
where 
\[\begin{split}
\tilde W_{\al,\be}(y)&=
\left\{\frac
{\GG\left(\frac{1}{2}+\be-\al\right)}
{\GG\left(\frac{1}{2}+\be+\al\right)}
\right\}^{1/2}W_{\al,\be}(y),\\\\
W_{\al,\be}(y)
&=\frac{e^{y/2}}{2\pi i}\int_{(\si)}\frac{\GG(w-\be)\GG(w+\be)}{\GG\left(\frac{1}{2}+w-\al\right)}\,y^{\frac{1}{2}-w}\,dw,\quad\si>|\Re\be|.
\end{split}\]
is the (normalized) Whittaker function. The normalization is introduced in order to retain the coefficients $\rho_\phi(\pm n)$ after the Maass operators have been applied. More precisely, if $k$ is an integer of the same parity as $\ka$, then
\[\phi_k(x+iy)=\sum_{n\neq 0}
\rfi(n)\tilde W_{\frac{n}{|n|}\frac{k}{2},i\mu}\bigl(4\pi|n|y\bigr)e(nx)\]
is a Maass form of weight $k$ and the same Petersson norm as $\phi$: 
\[\langle\phi_k,\phi_k\rangle=\langle\phi,\phi\rangle.\]
See Section~4 of \cite{DFI3} for details.

The unfolding technique yields an identity
\[(2\pi h)^{s-1}\int\limits_{\GG\backslash\HH}\phi_k(az)\ov\phi_k(bz)P_h(z,s)\frac{\,dx\,dy}{y^2}=\sum_{am-bn=h}\rfi(m)\brfi(n)H_{s,k,i\mu}\left(\frac{am+bn}{h}\right),\]
where
\[H_{s,k,i\mu}(u)=
\int_0^\infty\tilde W_{\frac{u+1}{|u+1|}\frac{k}{2},i\mu}\bigl(|u+1|y\bigr)\bar{\tilde W}_{\frac{u-1}{|u-1|}\frac{k}{2},i\mu}\bigl(|u-1|y\bigr)y^{s-2}\,dy,\quad u\neq\pm 1.\]

The main question that arises in the light of the above discussion is which weight functions $W:\RR\to\CC$ can be obtained by an averaging device from the $H_{s,k,i\mu}$ corresponding to values $s$ on a vertical line $\si+i\RR$ ($\si>1$) and all even (resp. odd) integers $k$. In Chapter~\ref{convolution} we shall make the first step in answering these questions by obtaining a fairly precise description of the span of the functions $H_{s,0,i\mu}$ (Theorem~\ref{Th91}).
\chapter{Approximate functional equation}\label{functional}

\section{Overview}

We shall approximate the values of a principal $L$-function $L(s,\pi)$ 
on the critical line $\Re s=\frac{1}{2}$ as a sum of 
two truncated Dirichlet series which have about $\sqrt{C(s,\pi)}$ terms.
We borrow notation from Section~\ref{size}, and we also refer the
reader to Section~\ref{functional_intro} for an introduction. The results of this chapter were published in \cite{Ha1}.

In order to keep the argument as clean as possible, we shall only display
our formulae for the central value $L\left(\frac{1}{2},\pi\right)$. This results in
no loss of generality, as $L\left(\frac{1}{2}+it,\pi\right)$ can be interpreted as
the central value corresponding to the twisted 
representation $\pi\otimes|\det|^{it}$. For convenient reference we record
the change of parameters in the formulae as we twist $\pi$ by a 1-dimensional
representation.
\[\pi\nyil\pi\otimes|\det|^{it};\quad
L\left(\tfrac{1}{2},\pi\right)\nyil L\left(\tfrac{1}{2}+it,\pi\right);\quad
C\left(\tfrac{1}{2},\pi\right)\nyil C\left(\tfrac{1}{2}+it,\pi\right);\]
\[\la_\pi(n)\nyil n^{-it}\la_\pi(n);\qquad\mu_j\nyil\mu_j-it;\qquad N\nyil N;
\qquad\ka\nyil N^{-it}\ka.\]

For the rest of this chapter $\pi$ will be a fixed cusp form on $\GL_m$ over
a number field $F$, and $C$ will abbreviate 
\begin{equation}\label{eq031}
C=C\left(\frac{1}{2},\pi\right)=
\frac{N}{(2\pi)^{md}}\prod_{j=1}^{md}\left|\frac{1}{2}-\mu_j\right|.
\end{equation}

\begin{theorem}\label{theorem01} There is a smooth function $f:(0,\infty)\to\CC$ 
and a complex number $\la$ of modulus 1
depending only on the Archimedean
parameters $\mu_j$ $(j=1,\dots,md)$ such that
\begin{equation}\label{eq08}
L\left(\frac{1}{2},\pi\right)=
\sum_{n=1}^\infty\frac{\la_\pi(n)}{\sqrt{n}}f\left(\frac{n}{\sqrt{C}}\right)+
\ka\la\sum_{n=1}^\infty\frac{\ov\la_\pi(n)}{\sqrt{n}}
\ov{f}\left(\frac{n}{\sqrt{C}}\right).
\end{equation}
The function $f$ and its partial derivatives $f^{(k)}$
$(k=1,2,.\dots)$ satisfy the following uniform growth estimates at
$0$ and infinity:
\begin{equation}\label{eq05}f(x)=
\begin{cases}1+O_{\si}(x^{\si}),&\quad
0<\si<\tfrac{1}{m^2+1};\\
O_{\si}(x^{-\si}),&\quad \si>0;
\end{cases}
\end{equation}
\begin{equation}\label{eq012}
f^{(k)}(x)=O_{\si,k}(x^{-\si}),\quad \si>k-\tfrac{1}{m^2+1}.
\end{equation}
The implied constants depend only on $\si$, $k$, $m$ and $d$.
\end{theorem}

\begin{Remark}
The range $0<\si<\frac{1}{m^2+1}$ in (\ref{eq05}) can be widened
to $0<\si<\frac{1}{2}$ for all representations $\pi$ which are tempered at
$\infty$, that is, conjecturally for all $\pi$. Similarly, upon the
Ramanujan--Selberg conjecture the range of $\si$ in (\ref{eq012})
can be extended to $\si>k-\frac{1}{2}$.
\end{Remark}

Combining the theorem with Molteni's bound (\ref{eq74}) we obtain that the
size of the central value $L\left(\frac{1}{2},\pi\right)$ can be very well
approximated with the first $C^{1/2+\ep}$ Dirichlet coefficients.

\begin{Cor}\label{Cor1} For any positive numbers $\ep$ and $A$,
\[L\left(\frac{1}{2},\pi\right)=
\sum_{n\leq
C^{1/2+\ep}}\frac{\la_\pi(n)}{\sqrt{n}}f\left(\frac{n}{\sqrt{C}}\right)+
\ka\la\sum_{n\leq C^{1/2+\ep}}\frac{\ov\la_\pi(n)}{\sqrt{n}}
\ov{f}\left(\frac{n}{\sqrt{C}}\right)+O_{\ep,A}(C^{-A}).\]
The implied constant depends only on $\ep$, $A$, $m$ and $d$.
\end{Cor}

In particular, by applying (\ref{eq74}) again, we can reconstruct the convexity bound 
(\ref{eq80}) for the central value (in fact for all values on the critical line). 

In a family of representations $\pi$, it is often desirable to see
that the weight functions $f$ do not vary too much. In fact,
assuming that the Archimedean parameters are not too small, one can
replace $f$ by an explicit function $g$ (independent of $\pi$) and
derive an approximate functional equation with a nontrivial error
term, that is, an error substantially smaller than the convexity
bound furnished by the above corollary. To state the result, we
introduce
\begin{equation}\label{eq019}
\eta=\min_{j=1,\dots,md}\left|\frac{1}{2}-\mu_j\right|.
\end{equation}

\begin{theorem}\label{theorem02} Let $g:(0,\infty)\to\RR$ be a smooth function
with the functional equation $g(x)+g(1/x)=1$ and derivatives decaying
faster than any negative power of $x$ as $x\to\infty$. Then, for
any $\ep>0$,
\[L\left(\frac{1}{2},\pi\right)=
\sum_{n=1}^\infty\frac{\la_\pi(n)}{\sqrt{n}}g\left(\frac{n}{\sqrt{C}}\right)+
\ka\la\sum_{n=1}^\infty\frac{\ov\la_\pi(n)}{\sqrt{n}}
g\left(\frac{n}{\sqrt{C}}\right)+O_{\ep,g}(\eta^{-1}C^{1/4+\ep}),\]
where $\la$ (of modulus 1) is given by (\ref{eq030}), and
the implied constant depends only on $\ep$, $g$, $m$ and $d$.
\end{theorem}

\begin{Remark} The formula is really of value when the family
under consideration satisfies $\eta\gg C^{\de}$ with some fixed
$\de>0$.
\end{Remark}

\section{The implicit form}

In this section we prove Theorem~\ref{theorem01}.
We introduce the auxiliary function
\begin{equation}\label{eq011}
F(s,\pi_\infty)=\left\{N^s
\frac{L\left(\frac{1}{2}+s,\pi_\infty\right)L\left(\frac{1}{2},\tilde\pi_\infty\right)}
{L\left(\frac{1}{2}-s,\tilde\pi_\infty\right)L\left(\frac{1}{2},\pi_\infty\right)}
\right\}^{1/2},
\end{equation}
which is holomorphic in the half plane $\Re s>-\frac{1}{m^2+1}$ by (\ref{eq81})
and (\ref{eq70}). With this notation
we can rewrite the functional equation (\ref{eq71}) as
\begin{equation}\label{eq01}
F(s,\pi_\infty)L\left(\tfrac{1}{2}+s,\pi\right)=\ka\la
F(-s,\tilde\pi_\infty)L\left(\tfrac{1}{2}-s,\tilde\pi\right),
\end{equation}
where
\begin{equation}\label{eq030}
\la=\frac{L\left(\frac{1}{2},\tilde\pi_\infty\right)}{L\left(\frac{1}{2},\pi_\infty\right)}.
\end{equation}
It follows from (\ref{eq82}) that $|\la|=1$, $F(0,\pi_\infty)=1$, and
\begin{equation}\label{eq06}
\ov{F}(s,\pi_\infty)=F(\ov{s},\tilde\pi_\infty).
\end{equation}

We also fix an entire function $H(s)$ which satisfies
the growth estimate
\begin{equation}\label{eq020}
H(s)\ll_{\si,A}\bigl(1+|s|\bigr)^{-A},\quad\Re s=\si;
\end{equation}
on vertical lines. In addition, we shall assume that $H(0)=1$ and that
$H(s)$ is symmetric with respect to both axes:
\begin{equation}\label{eq09}
H(s)=H(-s)=\ov{H}(\ov{s}).
\end{equation}
Such a function can be obtained as the Mellin transform of a
smooth function $h:~(0,\infty)\to\RR$ which has total mass 1 with
respect to the measure $dx/x$, functional equation $h(1/x)=h(x)$,
and derivatives decaying faster than any negative power of $x$ as
$x\to\infty$:
\[H(s)=\int_0^\infty h(x)x^s\frac{dx}{x}.\]

Using these two auxiliary functions and taking an arbitrary
$0<\si<\frac{1}{m^2+1}$, we can express the central value $L\left(\frac{1}{2},\pi\right)$
via the residue theorem as
\[\begin{split}
L\left(\frac{1}{2},\pi\right)
&=\frac{1}{2\pi i}\int_{(\si)}L\left(\frac{1}{2}+s,\pi\right)F(s,\pi_\infty)H(s)\frac{ds}{s}\\
&-\frac{1}{2\pi i}\int_{(-\si)}L\left(\frac{1}{2}+s,\pi\right)F(s,\pi_\infty)H(s)\frac{ds}{s}.
\end{split}\]
This step is justified by the convexity bound (\ref{eq80}), inequality 
(\ref{eq020}) and Lemma~\ref{Lemma2} below. Applying a
change of variable $s\mapsto -s$ in the second integral we get, by
the functional equations (\ref{eq01}) and (\ref{eq09}),
\[\begin{split}
L\left(\frac{1}{2},\pi\right)
&=\frac{1}{2\pi i}\int_{(\si)}L\left(\frac{1}{2}+s,\pi\right)F(s,\pi_\infty)H(s)\frac{ds}{s}\\
&+\frac{\ka\la}{2\pi i}
\int_{(\si)}L\left(\frac{1}{2}+s,\tilde\pi\right)F(s,\tilde\pi_\infty)H(s)\frac{ds}{s}.
\end{split}\]
The second integral is minus the complex conjugate of the first one, as can be seen
by another change of variable $s\mapsto\ov{s}$ combined with the functional equations
(\ref{eq82}), (\ref{eq06}) and (\ref{eq09}). Therefore we obtain the representation
(\ref{eq08}) of Theorem~\ref{theorem01} by defining
\begin{equation}\label{eq025}
f\left(\frac{x}{\sqrt{C}}\right)= \frac{1}{2\pi
i}\int_{(\si)}x^{-s}F(s,\pi_\infty)H(s)\frac{ds}{s}.
\end{equation}
For any nonnegative integer $k$ we also have
\begin{equation}\label{eq013}
f^{(k)}(x)=\frac{(-1)^k}{2\pi i}
\int_{(\si)}x^{-s-k}C^{-s/2}F(s,\pi_\infty)H(s)s(s+1)\dots(s+k-1)
\frac{ds}{s}.
\end{equation}

When $k=0$, the integrand in this expression is holomorphic for
$\Re s>-\frac{1}{m^2+1}$ with the exception of a simple pole at $s=0$
with residue 1. So in this case we are free to move the line of
integration to any nonzero $\si>-\frac{1}{m^2+1}$, but negative $\si$'s
will pick up an additional value 1 from the pole at $s=0$. When
$k>0$, the integrand is holomorphic in the entire half plane $\Re
s>-\frac{1}{m^2+1}$, so the line of integration can be shifted to any
$\si>-\frac{1}{m^2+1}$ without changing the value of the integral.
Henceforth, by (\ref{eq020}) and (\ref{eq013}), the truth of
inequalities (\ref{eq05}) and (\ref{eq012}) is reduced to the
following:

\begin{Lemma}\label{Lemma2} For any $\si>-\frac{1}{m^2+1}$, there is a uniform bound
\begin{equation}\label{eq022}
C^{-s/2}F(s,\pi_\infty)\ll_{\si}\bigl(1+|s|\bigr)^{md\si/2},\quad\Re s=\si.
\end{equation}
The implied constant depends only on $\si$, $m$ and $d$.
\end{Lemma}

We start with the following simple estimate.

\begin{Lemma}\label{Lemma1}
For any $\al>-\si$, there is a uniform bound
\[\frac{\GG(z+\si)}{\GG(z)}\ll_{\al,\si}|z+\si|^\si,\quad\Re z\geq\al.\]
\end{Lemma}
\noindent \emph{Proof of Lemma~\ref{Lemma1}.} The function
$\GG(z+\si)/\GG(z)$ is holomorphic in a neighborhood of $\Re
z\geq\al$. For $|z|>2|\si|$ we get, using Stirling's formula,
\[\frac{\GG(z+\si)}{\GG(z)}
\ll_\si\left|\frac{(z+\si)^{z+\si-1/2}}{z^{z-1/2}}\right|
\ll_\si|z+\si|^\si.\] The rest of the values of $z$ (those
with $\Re z\geq\al$ and $|z|\leq 2|\si|$) form a compact set, so
for these we simply have
\[\frac{\GG(z+\si)}{\GG(z)}\ll_{\al,\si}1\ll_{\al,\si}|z+\si|^\si.\qed\]

\noindent \emph{Proof of Lemma~\ref{Lemma2}.} 
Let $s=\si+it$. For any
$j=1,\dots,md$, we apply Lemma~\ref{Lemma1} with
\[\al=\frac{1}{2(m^2+1)}-\frac{\si}{2},\quad
z=\frac{1}{4}-\frac{\mu_j}{2}-\frac{\si}{2}+\frac{it}{2}\]
to see that
\[\frac
{\GG\left(\frac{1}{4}-\frac{\mu_j}{2}+\frac{\si}{2}+\frac{it}{2}\right)}
{\GG\left(\frac{1}{4}-\frac{\mu_j}{2}-\frac{\si}{2}+\frac{it}{2}\right)}
\ll_{\si,m}
\left|\frac{1}{4}-\frac{\mu_j}{2}+\frac{\si}{2}+\frac{it}{2}\right|^\si.\]
This is the same as
\[\frac
{\GG\left(\frac{1}{4}-\frac{\mu_j}{2}+\frac{s}{2}\right)}
{\GG\left(\frac{1}{4}-\frac{\ov\mu_j}{2}-\frac{s}{2}\right)}
\ll_{\si,m}
\left|\frac{1}{2}-\mu_j+s\right|^\si.\] 
It follows from (\ref{eq70}) that
\[\left|\frac{1}{2}-\mu_j+s\right|\leq
\left|\frac{1}{2}-\mu_j\right|+|s|\ll_m
\left|\frac{1}{2}-\mu_j\right|\bigl(1+|s|\bigr),\]
therefore we have
\[\frac
{\GG\left(\frac{1}{4}-\frac{\mu_j}{2}+\frac{s}{2}\right)}
{\GG\left(\frac{1}{4}-\frac{\ov\mu_j}{2}-\frac{s}{2}\right)}
\ll_{\si,m}\left|\frac{1}{2}-\mu_j\right|^\si\bigl(1+|s|\bigr)^\si.\]
Taking the product of these inequalities over all $j=1,\dots,md$,
and using (\ref{eq81}), (\ref{eq82}) and (\ref{eq031}), we get
\[\frac{L\left(\frac{1}{2}+s,\pi_\infty\right)}
{L\left(\frac{1}{2}-s,\tilde\pi_\infty\right)}
\ll_{\si,m,d}\left(\frac{C}{N}\right)^\si\bigl(1+|s|\bigr)^{md\si},\quad\Re s=\si.\] 
By (\ref{eq011}), this is equivalent to (\ref{eq022}), completing the
proof of Lemma~\ref{Lemma2} and Theorem~\ref{theorem01}.\qed

\section{The explicit form}

Our aim is to deduce Theorem~\ref{theorem02}.
We can assume that $H(s)$ is the Mellin transform of
$h(x)=-xg'(x)$. Indeed, $h:(0,\infty)\to\RR$ is a smooth function
with the functional equation $h(1/x)=h(x)$ and 
derivatives decaying
faster than any negative power of $x$ as $x\to\infty$, therefore
$H(s)$ is entire and satisfies (\ref{eq020}) and (\ref{eq09}). Also,
\[H(0)=-\int_0^\infty g'(x)=g(0+)=1.\]
Equivalently, $H(s)/s$ is the Mellin transform of $g(x)$, because
by partial integration it follows that
\[-\int_0^\infty g'(x)x^sdx=s\int_0^\infty g(x)x^s\frac{dx}{x}.\]
In any case, $g(x)$ can be expressed as an inverse Mellin transform
\[g(x)=\frac{1}{2\pi i}\int_{(\si)}x^{-s}H(s)\frac{ds}{s}.\]

The idea is to compare $g(x)$ with the function $f(x)$ given by
(\ref{eq025}). We have, for any $\si>0$,
\[f(x)-g(x)=\frac{1}{2\pi i}\int_{(\si)}x^{-s}
\bigl\{C^{-s/2}F(s,\pi_\infty)-1\bigr\}H(s)\frac{ds}{s}.\]
In fact, the integrand is holomorphic in
the entire half plane $\Re s>-\frac{1}{m^2+1}$,
so the line of integration can be shifted to any
$\si>-\frac{1}{m^2+1}$ without changing the value of the integral.
In particular, the choice $\si=0$ is permissible, that is,
\begin{equation}\label{eq024}
f(x)-g(x)=\frac{1}{2\pi i}\int_{-\infty}^\infty x^{-it}
\bigl\{C^{-it/2}F(it,\pi_\infty)-1\bigr\}H(it)\frac{dt}{t}.
\end{equation}

Note that $x^{-it}$ and $C^{-it/2}F(it,\pi_\infty)$ are of modulus
1. For any $\ep>0$, the values of $t$ with
$|t|\geq\min(\eta/2,C^\ep)$ contribute $O_{\ep,g,m,d}(\eta^{-1})$ to
the integral. This follows from (\ref{eq020}) and $\eta\ll
C^{1/md}$. We estimate the remaining contribution via the following lemma.

\begin{Lemma}\label{Lemma3} For any $\ep>0$, there is a uniform bound
\[C^{-it/2}F(it,\pi_\infty)-1\ll_\ep|t|\eta^{-1}C^\ep
,\quad |t|<\min(\eta/2,C^\ep).\] The implied constant depends only
on $\ep$, $m$ and $d$.
\end{Lemma}

\noindent \emph{Proof.}
As $C^{-it/2}F(it,\pi_\infty)$ lies on the unit circle,
it suffices to show that
\[\log\bigl\{C^{-it/2}F(it,\pi_\infty)\bigr\}
\ll_{\ep,m,d}|t|\eta^{-1}C^\ep ,\quad |t|<\min(\eta/2,C^\ep).\] Here
the left hand side is understood as a continuous function defined
via the principal branch of the logarithm near $t=0$. Using
 (\ref{eq031}), (\ref{eq011}), (\ref{eq81}) and (\ref{eq82}) we can see that the
derivative (with respect to $t$) of the left hand side is given by
\[\frac{i}{2}\Re\sum_{j=1}^{md}\left\{\frac{\GG'}{\GG}
\left(\frac{1}{4}-\frac{\mu_j}{2}+\frac{it}{2}\right)
-\log\left(\frac{1}{4}-\frac{\mu_j}{2}\right)\right\},\] so we can
further reduce the lemma to
\begin{equation}\label{eq018}
\frac{\GG'}{\GG}\left(\frac{1}{4}-\frac{\mu_j}{2}+\frac{it}{2}\right)
-\log\left(\frac{1}{4}-\frac{\mu_j}{2}\right)\ll_{\ep,m,d}\eta^{-1}C^\ep
,\quad |t|<\min(\eta/2,C^\ep).
\end{equation}
Here $\frac{1}{4}-\frac{\mu_j}{2}+\frac{it}{2}$ has real part at least 
$\frac{1}{2(m^2+1)}$ by (\ref{eq70})
and absolute value at least $\eta/4$ by (\ref{eq019}).
Therefore, a standard bound yields
\[\frac{\GG'}{\GG}\left(\frac{1}{4}+\frac{\mu_j}{2}+\frac{it}{2}\right)=
\log\left(\frac{1}{4}+\frac{\mu_j}{2}+\frac{it}{2}\right)+O_m(\eta^{-1}).\]
For $|t|<\min(\eta/2,C^\ep)$ we can also see that
\[\log\left(\frac{1}{4}+\frac{\mu_j}{2}+\frac{it}{2}\right)=
\log\left(\frac{1}{4}+\frac{\mu_j}{2}\right)+O(\eta^{-1}C^\ep).\]
It follows from (\ref{eq70}) that $C\gg_{m,d}1$, therefore the
last two estimates add up to (\ref{eq018}) as required.\qed

\medskip

Returning to the integral (\ref{eq024}), it follows from
Lemma~\ref{Lemma3} that the values of $t$ with
$|t|<\min(\eta/2,C^\ep)$ contribute at most
$O_{\ep,g,m,d}(\eta^{-1}C^{2\ep})$. Altogether we have, by
$C\gg_{m,d}1$,
\[f(x)-g(x)=O_{\ep,g,m,d}(\eta^{-1}C^{2\ep}).\]

We conclude Theorem~\ref{theorem02} by combining this estimate with
Corollary~\ref{Cor1} and Molteni's bound (\ref{eq74}).

\chapter{Shifted convolution sums and the circle method}\label{maass}

\section{Overview}\label{sect6}

We shall establish, in the spirit of Duke, Friedlander and Iwaniec, a nontrivial bound for the shifted convolution sums (\ref{eq100}) arising from classical holomorphic or Maass cusp forms for the Hecke congruence subgroups.
We refer the reader to Section~\ref{shifted} for an introduction. The notions in the following theorem will be defined in the next section. The result, in less explicit form, will also appear in \cite{Ha2}.

\begin{theorem}\label{Th1} Let $\lfi(m)$ (resp. $\lpsi(n)$)
be the normalized Fourier coefficients of a holomorphic or Maass
cusp form $\phi$ (resp. $\psi$) of level $N$ and arbitrary nebentypus character modulo $N$. Let $|\tilde\mu|$ (resp. $|\tilde\nu|$) denote the Archimedean size of $\phi$ (resp. $\psi$), and suppose that $f$ satisfies (\ref{eq1}). Then for coprime $a$ and $b$ we have
\[D_f(a,b;h)\ll P^{11/10}N^{9/5}|\tilde\mu\tilde\nu|^{9/5+\ep}(ab)^{-1/10}(X+Y)^{1/10}(XY)^{2/5+\ep},\]
where the implied constant depends only on $\ep$.
\end{theorem}

\begin{Remark} We shall see in Section~\ref{sect2} that Cauchy's
inequality implies
\begin{equation}\label{eq40}
D_f(a,b;h)\ll N|\tilde\mu\tilde\nu|^{1/2}(ab)^{-1/2}(XY)^{1/2}.
\end{equation}
The conclusion of the theorem supercedes this trivial
bound whenever
\begin{equation}\label{eq17}
P^{11}N^8|\tilde\mu\tilde\nu|^{13+\ep}(ab)^4\ll\frac{(XY)^{1-\ep}}{X+Y}.
\end{equation}
\end{Remark}

The proof of Theorem~\ref{Th1} is presented in Sections~\ref{sect1}
through \ref{bessel} and closely follows \cite{DFI}. The heart of
the argument is a Vorono\"\i-type summation formula (see
Section~\ref{sect1b}) for transforming certain exponential sums
defined by the coefficients $\lfi(m)$ and $\lpsi(n)$. As
the level of the forms imposes some restriction on the frequencies
in the formula, we replace (in Section~\ref{sect3})
the classical Farey dissection (or the
$\dd$-method) with Jutila's variant of the circle method \cite{Ju3}.
The variant uses overlapping intervals, and hence provides great
flexibility in the choice of frequencies.
After transforming our exponential generating functions
in Section~\ref{sect10}, we encounter twisted Kloosterman sums
\[S_\chi(m,n;q)=\csillag\sum_{d\smod{q}}
\chi(d)e_q\bigl(dm+{\bar d}n\bigr),\] where $\chi$ is a Dirichlet
character mod $q$. We refer to the usual Weil--Estermann
bound \begin{equation}\label{eq10}\bigl|S_\chi(m,n;q)\bigr|\leq
(m,n,q)^{1/2}q^{1/2}\tau(q),\end{equation}
for which the original proofs \cite{W,E} can be adapted.
In Section~\ref{sect2} we apply a smooth dyadic decomposition,
and conclude the theorem by optimizing
the free parameters. In order to achieve polynomial uniformity in the Archimedean parameters of the cusp forms, we need to exhibit careful estimates for the Bessel functions involved in the summation formula. These estimates appear in Section~\ref{bessel} with detailed proofs.

\section{Normalized Fourier coefficients}\label{sect1}

We define the normalized Fourier coefficients of cusp forms as
follows. Let $\phi$ be a cusp form of level $N$ and nebentypus
$\chi$, that is, a holomorphic cusp form of some integral weight
$k$, or a real-analytic Maass cusp form of some nonnegative Laplacian
eigenvalue $1/4+\mu^2$. 
In the holomorphic case we write $k-1=2i\mu$,
in the real-analytic case we define $k=0$, 
and in both cases we put $\tilde\mu=1/2+i\mu$ and call $|\tilde\mu|$ the 
\emph{Archimedean size} of $\phi$. This is in accordance with 
Section~\ref{size}.

By definition, $\chi$ is a Dirichlet
character mod $N$, and $\phi$ is a complex valued function on
the upper half plane $\HH=\{z:\Im z>0\}$, which decays exponentially to
zero at each cusp and satisfies a transformation
rule with respect to the Hecke congruence subgroup $\GG_0(N)$:
\[\phi\left(\frac{az+b}{cz+d}\right)=\chi(d)(cz+d)^{k}\phi(z),
\qquad\begin{pmatrix}a&b\\c&d\end{pmatrix}\in\GG_0(N).\]
In particular, $\phi$ admits the Fourier expansion
\begin{equation}\label{eq210}
\phi(x+iy)=\sum_{n\neq 0}{\rho}_\phi(n)W(ny)e(nx),\end{equation}
where
\begin{equation}\label{eq211}
W(y)=\begin{cases}
e^{-2\pi y}&\text{if $\phi$ is holomorphic,}\\
|y|^{1/2}K_{i\mu}\bigl(2\pi|y|\bigr)&\text{if
$\phi$ is real-analytic.}
\end{cases}\end{equation}
Here $e(x)=e^{2\pi ix}$, and $K_{i\mu}$ is the
MacDonald-Bessel function. If $\phi$ is holomorphic,
${\rho}_\phi(n)$ vanishes for $n<0$. Writing
\[\pet=\int_{\GG_0(N)\backslash\HH}y^{k-2}|\phi(x+iy)|^2\,dx\,dy,\]
we define the \emph{normalized Fourier coefficients} of $\phi$ as
\begin{equation}\label{eq200}
\lfi(n)=\begin{cases}
\left(\frac{N(k-1)!}{\pet(4\pi n)^{k-1}}\right)^{1/2}{\rho}_\phi(n)&\text{if $\phi$ is holomorphic,}\\
\left(\frac{N(4\pi|n|)}{\pet\cosh\pi\mu}\right)^{1/2}{\rho}_\phi(n)&\text{if $\phi$ is
real-analytic.}
\end{cases}\end{equation}
This normalization corresponds to Rankin--Selberg theory which implies the
following mean square estimate for the normalized Fourier coefficients (see Section~8.2 of \cite{I}):
\[c_N\sum_{1\leq|n|\leq x}|\lfi(n)|^2\sim x\qquad\text{as $x\to\infty$},\]
\[1\ll c_N\ll\log\log(3N).\]
More precisely,
\[c_N\asymp\frac{\vol\bigl(\GG_0(N)\backslash\HH\bigr)}{N}=
\frac{\pi}{3}\prod_{p\mid N}\left(1+\frac{1}{p}\right).\]
We also have a good uniform upper bound for all $x>0$ (see
Theorem~3.2 and (8.7) and (9.34) in \cite{I}): 
\begin{equation}\label{eq22}
\sum_{1\leq|n|\leq x}|\lfi(n)|^2\ll x+N|\tilde\mu|,
\end{equation}
where the implied constant is absolute.

\begin{Lemma}\label{bound}For any $\ep>0$ there is a uniform bound
\[y^{k/2}\phi(x+iy)\ll\pet^{1/2}\tmu^{3/2+\ep}y^{-\ep},
\quad x\in\RR,\ y>1/2.\]
The implied constant depends only on $\ep$.
\end{Lemma}

\noindent\emph{Proof.} We distungish between two cases.

\begin{Case2} $\phi$ is holomorphic. By (\ref{eq210}), (\ref{eq211}) and (\ref{eq200}), the statement is equivalent to
\[\left|\sum_{n=1}^\infty\lfi(n)(4\pi n)^\frac{k-1}{2}e^{-2\pi ny}e(nx)\right|^2\ll_\ep (k-1)!k^{3+\ep}Ny^{-k-\ep}.\]
By the Cauchy-Schwartz inequality the left hand side can be estimated from above by
\[\left(\sum_{n=1}^\infty|\lfi(n)|^2(4\pi n)^{-1-\ep}\right)
\left(\sum_{n=1}^\infty(4\pi n)^{k+\ep}e^{-4\pi ny}\right).\]
The first factor is $\ll_\ep Nk$ by the mean square bound (\ref{eq22}), therefore it remains to show that
\[\sum_{n=1}^\infty(4\pi ny)^{k+\ep}e^{-4\pi ny}\ll_\ep
(k-1)!k^{2+\ep}.\]
We accomplish this in stronger form by comparing the sum with the similar integral (note that $y\gg 1$):
\[\begin{split}
\sum_{n=1}^\infty(4\pi ny)^{k+\ep}e^{-4\pi ny}
&\ll\sup_{y>0}\left\{(4\pi ny)^{k+\ep}e^{-4\pi ny}\right\}+\int_0^\infty(4\pi ny)^{k+\ep}e^{-4\pi ny}\,dy\\\\
&=\left(\frac{k+\ep}{e}\right)^{k+\ep}+\GG(k+1+\ep)\ll_\ep (k-1)!k^{1+\ep}.
\end{split}\]
\end{Case2}

\begin{Case2}$\phi$ is real-analytic. By (\ref{eq210}), (\ref{eq211}) and (\ref{eq200}), the statement is equivalent to
\[\left|\sum_{n\neq 0}\lfi(n)K_{i\mu}\bigl(2\pi|n|y\bigr)e(nx)\right|^2\ll_\ep e^{-\pi|\mu|}\tmu^{3+\ep}Ny^{-1-\ep}.\]
By the Cauchy-Schwartz inequality the left hand side can be estimated from above by
\[\left(\sum_{n\neq 0}|\lfi(n)|^2|2\pi n|^{-1-\ep}\right)
\left(\sum_{n\neq 0}|2\pi n|^{1+\ep}\bigl|K_{i\mu}\bigl(2\pi|n|y\bigr)\bigr|^2\right).\]
The first factor is $\ll_\ep N\tmu$ by the mean square bound (\ref{eq22}), therefore it remains to show that
\[\sum_{n\neq 0}|2\pi ny|^{1+\ep}e^{\pi|\mu|}\bigl|K_{i\mu}\bigl(2\pi|n|y\bigr)\bigr|^2\ll_\ep
\tmu^{2+\ep}.\]
We accomplish this by employing Proposition~\ref{bessel2} of Section~\ref{bessel}, noting also that $y\gg 1$ and $|\Re(i\mu)|\leq 1/2$:
\begin{alignat*}{3}
\sum_{n\neq 0}|2\pi ny|^{1+\ep}e^{\pi|\mu|}\bigl|K_{i\mu}\bigl(2\pi|n|y\bigr)\bigr|^2
&=&
\sum_{4|n|y<\tmu}\dots\ \ 
&+&
\sum_{\tmu\leq 4|n|y<2\tmu}\dots\ 
&+
\sum_{2\tmu\leq 4|n|y}\dots\\\\
&\ll_\ep& 
\tmu^{2+\ep}\quad\ \ \,
&+&
\tmu^{1+\ep}\qquad\ 
&+
\quad\tmu^{\ep}.\qed
\end{alignat*}
\end{Case2}

\begin{Lemma}\label{bound2}For any $\ep>0$ there is a uniform bound
\[\bigl\|y^{k/2}\phi(x+iy)\bigr\|_\infty\ll\pet^{1/2}\tmu^{3/2+\ep}.\]
The implied constant depends only on $\ep$.
\end{Lemma}

\noindent\emph{Proof.} It is known that any $z=x+iy$ can be represented as
$z=\frac{aw+b}{cw+d}$, where 
$\left(\begin{smallmatrix}a&b\\c&d\end{smallmatrix}\right)\in\SL_2(\ZZ)$
and $w$ has imaginary part $\Im w>1/2$. The proof of the previous lemma can be adapted
almost verbatim 
to the cusp form $w\mapsto(cw+d)^{-k}\phi\left(\frac{aw+b}{cw+d}\right)$, so that
we have, in particular,
\[y^{k/2}\phi(x+iy)=\frac{|\Im w|^{k/2}}{|cw+d|^k}\phi\left(\frac{aw+b}{cw+d}\right)
\ll_\ep\pet^{1/2}\tmu^{3/2+\ep}.\qed\]

The proof of Theorem~\ref{Th1} is based on exponential sums of the form
\begin{equation}\label{eq213}
T_{\phi,\al}(M)=\sum_{1\leq m\leq M}\lfi(m)e(\al m).\end{equation}
We shall use the following uniform variant of Wilton's classical estimate.

\begin{proposition}\label{exponential} For any $\ep>0$ there is a uniform bound
\begin{equation}\label{eq212}
T_{\phi,\al}(M)\ll N^{1/2}\tmu^{2+\ep}M^{1/2+\ep},\quad\al\in\RR,\ M>0.
\end{equation}
The implied constant depends only on $\ep$.
\end{proposition}

\noindent\emph{Proof.} We can clearly assume that $M$ is a positive integer.
As before, we distungish between two cases.

\begin{Case3} $\phi$ is holomorphic. 
For any positive integer $m$ we have, by (\ref{eq210}), (\ref{eq211}) and (\ref{eq200}),
\[\lfi(m)e^{-2\pi my}e(m\al)=\left(\frac{N(k-1)!}{\pet(4\pi m)^{k-1}}\right)^{1/2}
\int_0^1\phi(\al+\be+iy)e(-m\be)\,d\be.\]
We multiply both sides by $(2\pi y)^{k/2+\ep}$,
and integrate with respect to $dy/y$. We obtain
\[\lfi(m)m^{-1/2-\ep}e(m\al)=\int_0^1\Phi_\al(\be)e(-m\be)\,d\be,\]
where
\[\Phi_\al(\be)=\frac{(2\pi)^{k/2+\ep}}{\GG(k/2+\ep)}\left(\frac{N(k-1)!}{\pet(4\pi)^{k-1}}\right)^{1/2}\int_0^\infty y^{k/2+\ep}\phi(\al+\be+iy)\frac{dy}{y}.\]
Note that the integral converges by Lemmata~\ref{bound} and \ref{bound2} and satisfies the uniform bound
\[\int_0^\infty y^{k/2+\ep}\phi(\al+\be+iy)\frac{dy}{y}\ll_\ep\pet^{1/2}k^{3/2+2\ep}.\]
It follows that
\[\Phi_\al(\be)\ll_\ep N^{1/2}k^{7/4+\ep}.\]
By introducing the kernel
\[F_M(\be)=\sum_{|m|\leq M}e(m\be)=\frac{\sin\pi(2M+1)\be}{\sin\pi\be},\]
we can write
\[\sum_{m=1}^M\lfi(m)m^{-1/2-\ep}e(m\al)=\int_0^1\Phi_\al(\be)F_M(\be)\,d\be.\]
It is known that the $\mathcal{L}^1$-norm of $F_M$ is $\ll\log(2M)$, therefore
it follows that
\[\sum_{m=1}^M\lfi(m)m^{-1/2-\ep}e(m\al)\ll_\ep N^{1/2}k^{7/4+\ep}M^\ep.\]
From this bound (\ref{eq212}) follows by partial summation.
\end{Case3}

\begin{Case3}$\phi$ is real-analytic.
For any nonzero integer $m$ we have, by (\ref{eq210}), (\ref{eq211}) and (\ref{eq200}),
\[\lfi(m)y^{1/2}K_{i\mu}(2\pi|m|y)e(m\al)=
\left(\frac{4\pi N}{\pet\cosh\pi\mu}\right)^{1/2}
\int_0^1\phi(\al+\be+iy)e(-m\be)\,d\be.\]
We multiply both sides by $(2\pi)^{1/2+\ep}y^\ep$,
and integrate with respect to $dy/y$. We obtain
\[\lfi(m)|m|^{-1/2-\ep}e(m\al)=\int_0^1\Phi_\al(\be)e(-m\be)\,d\be,\]
where
\[\Phi_\al(\be)=\frac{8\pi^{1+\ep}}{\prod\limits_\pm
\GG\left(\frac{1}{4}+\frac{\ep}{2}\pm\frac{i\mu}{2}\right)}
\left(\frac{N}{\pet\cosh\pi\mu}\right)^{1/2}
\int_0^\infty y^{\ep}\phi(\al+\be+iy)\frac{dy}{y}.\]
Note that the integral converges by Lemmata~\ref{bound} and \ref{bound2} and satisfies the uniform bound
\[\int_0^\infty y^{\ep}\phi(\al+\be+iy)\frac{dy}{y}\ll_\ep\pet^{1/2}\tmu^{3/2+2\ep}.\]
It follows that
\[\Phi_\al(\be)\ll_\ep N^{1/2}\tmu^{2+\ep},\]
and from this point we proceed exactly as in Case 1.\qed
\end{Case3}

\section{Summation formula}\label{sect1b}

Various Vorono\"\i-type summation formulas are fulfilled by the
normalized Fourier
coefficients. In the case of full level ($N=1$) Duke and Iwaniec
\cite{DI} established such a formula for holomorphic cusp forms
and Meurman \cite{M} for Maass cusp forms. These can be
generalized to arbitrary level and nebentypus with obvious minor
modifications as follows.
\begin{proposition}\label{Prop1}
Let $d$ and $q$ be coprime integers such that $N\mid q$, and let
$g$ be a smooth, compactly supported function on $(0,\infty)$. If
$\phi$ is a holomorphic cusp form of level $N$, nebentypus $\chi$
and integral weight $k$ then
\[\chi(d)\sum_{n=1}^\infty\lfi(n)e_q(dn)g(n)=
\sum_{n=1}^\infty\lfi(n)e_q\bigl(-{\bar d}n\bigr){\hat g}(n),\]
where
\[{\hat g}(y)=\frac{2\pi i^k}{q}\int_0^\infty g(x)
J_{k-1}\left(\frac{4\pi\sqrt{xy}}{q}\right)\,dx.\] If $\phi$ is a
real-analytic Maass cusp form of level $N$, nebentypus $\chi$ and
nonnegative Laplacian eigenvalue $1/4+\mu^2$ then
\[\chi(d)\sum_{n=1}^\infty\lfi(n)e_q(dn)g(n)=
\sum_{\pm}\sum_{n=1}^\infty\lfi(\mp n)e_q\bigl(\pm\bar d
n\bigr)g^{\pm}(n),\] where
\[\begin{split}
g^-(y)&=-\frac{\pi}{q\cosh\pi\mu}\int_0^\infty g(x)
\{Y_{2i\mu}+Y_{-2i\mu}\}\left(\frac{4\pi\sqrt{xy}}{q}\right)\,dx,\\\\
g^+(y)&=\frac{4\cosh\pi\mu}{q}\int_0^\infty g(x)
K_{2i\mu}\left(\frac{4\pi\sqrt{xy}}{q}\right)\,dx.
\end{split}\]
Here $\bar d$ is a multiplicative inverse of $d\bmod q$,
$e_q(x)=e(x/q)=e^{2\pi ix/q}$ and $J_{k-1}$, $Y_{\pm 2i\mu}$,
$K_{2i\mu}$ are Bessel functions.
\end{proposition}

The proof for the holomorphic case \cite{DI} is a straightforward
application of Laplace transforms. Meurman's proof for the
real-analytic case \cite{M} is more involved, but only because he
considers a wider class of test functions $g$ and has to deal with
delicate convergence issues. For smooth, compactly supported
functions $g$ as in our formulation these difficulties do not
arise, and one can give a much simpler proof based on Mellin
transformation, the functional equations of the $L$-series
attached to additive twists of $\phi$ (see \cite{M}), and Barnes'
formulae for the gamma function. Indeed, Lemma 5 in \cite{St}, a
special case of Meurman's summation formula, has been proved by
such an approach. We expressed the formula for the non-holomorphic
case in terms of $K$- and $Y$-Bessel functions in order to
emphasize the analogy with the Vorono\"\i-type formula for the
divisor function (where one has $\mu=0$) as derived by Jutila
\cite{Ju4,Ju5}. 

Michel recently extended the above formula to all denominators \cite{Mi,Mi2}. The extension becomes quite involved when $N$ is not square-free, and the proof relies heavily on Atkin--Lehner theory \cite{AL,Li,ALi}. We shall not use this general version.

\section{Setting up the circle method}\label{sect3}

For sake of exposition we shall only present the case of
Maass forms and the equation $am-bn=h$. The other
cases follow along similar lines by changing Bessel functions and
signs at relevant places of the argument. In our inequalities
$\ep$ will always denote a small positive number whose actual
value is allowed to change at each occurrence. Implied constants will 
always depend on $\ep$. All other dependencies will be explicitly indicated.

Let $\phi$ (resp. $\psi$) be a Maass cusp form of level $N$,
nebentypus $\chi$ (resp. $\om$) and Laplacian eigenvalue
$1/4+\mu^2\geq 0$ (resp. $1/4+\nu^2\geq 0$) whose normalized
Fourier coefficients are $\lfi(m)$ (resp. $\lpsi(n)$).
We shall first investigate $D_g(a,b;h)$ for smooth test functions
$g(x,y)$ which are supported in a box $[A,2A]\times[B,2B]$ and
have partial derivatives bounded by
\begin{equation}\label{eq27}g^{(k,l)}\ll_{k,l}A^{-k}B^{-l}P^{k+l}.
\end{equation} Our aim is to prove the estimate
\begin{equation}\label{eq29}
D_g(a,b;h)\ll P^{11/10}N^{9/5}|\tilde\mu\tilde\nu|^{9/5+\ep}(ab)^{-1/10}(A+B)^{1/10}(AB)^{2/5+\ep}.\end{equation}
In Section~\ref{sect2} we shall deduce Theorem~\ref{Th1} from this
bound by employing a partition of unity and decomposing
appropriately any smooth test function $f(x,y)$ satisfying
(\ref{eq1}). In fact, (\ref{eq29}) is a special case of
Theorem~\ref{Th1}, as can be seen upon setting $X=A$, $Y=B$, and
$f(x,y)=g(x,y)$. 

We shall assume that
\begin{equation}\label{eq28}
P^{11}N^8|\tilde\mu\tilde\nu|^{13+\ep}(ab)^4\ll\frac{(AB)^{1-\ep}}{A+B},
\end{equation}
for otherwise (\ref{eq29}) follows from the trivial upper bound
\begin{equation}\label{eq107}
D_g(a,b;h)\ll N|\tilde\mu\tilde\nu|^{1/2}(ab)^{-1/2}(AB)^{1/2}.
\end{equation}
The trivial bound itself is a consequence of $g\ll 1$, Cauchy's
inequality, and the mean square estimate (\ref{eq22}) applied
to the forms $\phi$ and $\psi$.

As $g(x,y)$ is supported in $[A,2A]\times[B,2B]$, we can assume
that $A,B\geq 1/2$, and also that
\begin{equation}\label{eq11}
h\leq 2(A+B),
\end{equation}
for otherwise $D_g(a,b;h)$ vanishes trivially. We shall attach, as
in \cite{DFI}, a redundant factor $w(x-y-h)$ to $g(x,y)$, where
$w(t)$ is a smooth function supported on $|t|\leq\dd^{-1}$ such that
$w(0)=1$ and $w^{(i)}\ll_i\dd^i$. This, of course, does not alter
$D_g(a,b;h)$. We choose
\begin{equation}\label{eq2}
\dd=P\frac{A+B}{AB},
\end{equation}
so that, by (\ref{eq27}), the new function
\[F(x,y)=g(x,y)w(x-y-h)\]
satisfies
\begin{equation}\label{eq102}
|x-y-h|>\dd^{-1}\quad\Longrightarrow\quad F(x,y)=0,
\end{equation}
and its partial derivatives are bounded by
\begin{equation}\label{eq5}
F^{(k,l)}\ll_{k,l}\dd^{k+l}.
\end{equation}
We apply the Hardy--Littlewood method to detect the equation $am-bn=h$,
that is, we express $D_F(a,b;h)$ as the integral of a certain
exponential sum over the unit interval $[0,1]$. We get
\begin{equation}\label{eq30}D_g(a,b;h)=D_F(a,b;h)=\int_0^1
G(\al)\,d\al,\end{equation} where
\[G(\al)=\sum_{m,n}\lfi(m)\lpsi(n)F(am,bn)e\bigl((am-bn-h)\al\bigr).\]
We shall approximate this integral by the following proposition of Jutila
(a consequence of the main theorem in \cite{Ju3}).
\begin{proposition}[Jutila]\label{Prop2}
Let $\CQ$ be a nonempty set of integers $Q\leq q\leq 2Q$, where $Q\geq 1$.
Let $Q^{-2}\leq\dd\leq Q^{-1}$, and for
each fraction $d/q$ (in its lowest terms) denote by
$I_{d/q}(\al)$ the characteristic function of the interval
$\left[d/q-\dd,d/q+\dd\right]$. Write $L$ for the
number of such intervals, that is,
\[L=\sum_{q\in\CQ}\varphi(q),\]
and put
\[\tilde{I}(\al)=\frac{1}{2\dd L}\sum_{q\in\CQ}\ \
\csillag\sum_{d\smod{q}}I_{d/q}(\al).\] If $I(\al)$ is the
characteristic function of the unit interval $[0,1]$, then
\[\int_{-\infty}^\infty\bigl(I(\al)-\tilde{I}(\al)\bigr)^2 \,dx\ll
\dd^{-1}L^{-2}Q^{2+\ep},\]
where the implied constant depends on $\ep$ only.
\end{proposition}

We shall choose some $Q$ and apply the proposition with a set of
denominators of the form
\[\CQ=\bigl\{q\in[Q,2Q]:Nab\mid q\text{ and
}(h,q)=(h,Nab)\bigr\}.\] 
By a result of Jacobsthal \cite{Ja}, the
largest gap between reduced residue classes mod $h$ is of size
$\ll h^\ep$, whence (\ref{eq11}) shows that
\begin{equation}\label{eq18}
|\CQ|\gg\frac{Q(AB)^{-\ep}}{Nab},
\end{equation}
assuming the right hand side exceeds some large positive constant
$c=c(\ep)\geq 1$. Moreover, we shall assume that
\begin{equation}\label{eq3}
Q^{-2}\leq\dd\leq Q^{-1},
\end{equation}
so that also
\begin{equation}\label{eq12}
1\leq Q\leq AB,
\end{equation}
whence (\ref{eq18}) yields
\begin{equation}\label{eq6}
L\gg\frac{Q^2(AB)^{-\ep}}{Nab}.
\end{equation}
We clearly have
\begin{equation}\label{eq4}
|D_F(a,b;h)-\tilde{D}_F(a,b;h)|\leq\|G\|_\infty\|I-\tilde{I}\|_1,
\end{equation}
where
\[\begin{split}
\tilde{D}_F(a,b;h)&=\int_{-\infty}^\infty
G(\al)\tilde{I}(\al)\,d\al =\frac{1}{2\dd L}\sum_{q\in\CQ}\ \
\csillag\sum_{d\smod{q}}
\int_{-\infty}^\infty G(\al)I_{d/q}(\al)\,d\al\\\\
&=\frac{1}{2\dd L}\sum_{q\in\CQ}\ \ \csillag\sum_{d\smod{q}}
\int_{-\dd}^\dd G(d/q+\beta)\,d\beta =\frac{1}{2\dd
L}\sum_{q\in\CQ}\ \ \csillag\sum_{d\smod{q}}\FI_{d/q},
\end{split}\]
say. To derive an upper estimate for $G(\al)$, we express it as
\[G(\al)=\int_0^\infty\int_0^\infty F(x,y)e(-h\al)
\,dT_{\phi,a\al}(x/a)\,dT_{\psi,-b\al}(y/b),\]
where the exponential sums $T_{\phi,a\al}$ and $T_{\psi,-b\al}$ are defined by
(\ref{eq213}).
Using the uniform bound provided by Proposition~\ref{exponential}, and
also (\ref{eq102}) and (\ref{eq5}), it follows that
\[\|G\|_\infty\ll\frac{N|\tilde\mu\tilde\nu|^{2+\ep}}{(ab)^{1/2}}(AB)^{1/2+\ep}\|F^{(1,1)}\|_1
\ll\frac{N|\tilde\mu\tilde\nu|^{2+\ep}\dd}{(ab)^{1/2}}\cdot\frac{(AB)^{3/2+\ep}}{A+B}.\] Also, by
(\ref{eq6}) and Proposition~\ref{Prop2} we get
\[\|I-\tilde{I}\|_1\leq 3\|I-\tilde{I}\|_2\ll \frac{Nab}{\dd^{1/2}Q}(AB)^\ep,\]
so that (\ref{eq4}) becomes
\begin{equation}\label{eq15}
D_F(a,b;h)-\tilde{D}_F(a,b;h)\ll 
\frac{N^{2}|\tilde\mu\tilde\nu|^{2+\ep}(ab)^{1/2}\dd^{1/2}}{Q}\cdot\frac{(AB)^{3/2+\ep}}{A+B}.
\end{equation}

\section{Transforming exponential sums}\label{sect10}

The contribution of the interval $[d/q-\dd,d/q+\dd]$ can be
expressed as
\[\FI_{d/q}=\int_{-\dd}^\dd
G(d/q+\beta)\,d\beta=e_q(-dh)\sum_{m,n}\lfi(m)\lpsi(n)e_q\bigl(d(am-bn)\bigr)E(m,n),\]
where
\begin{equation}\label{eq109}
E(x,y)=F(ax,by)\int_{-\dd}^\dd
e\bigl((ax-by-h)\beta\bigr)\,d\beta.
\end{equation}
For further reference we record the following two simple consequences of
(\ref{eq102}) and (\ref{eq5}):
\[E^{(k,l)}\ll_{k,l}\dd^{k+l+1}a^kb^l;\]
\begin{equation}\label{eq7}
\|E^{(k,l)}\|_1\ll_{k,l}\dd^{k+l}a^{k-1}b^{l-1}\frac{AB}{A+B}.
\end{equation}
We assume that $q\in\CQ$, hence $Nab\mid q$, and Proposition~\ref{Prop1} 
yields
\[\FI_{d/q}=\overline{\chi\om}(d)e_q(-dh)\sum_{\pm\pm}\,\sum_{m,n\geq
1}\lfi(\mp m)\lpsi(\mp n)e_q\bigl({\bar d}(\pm am\mp
bn)\bigr)E^{\pm\pm}(m,n),\] where the corresponding signs must be
matched, and
\[E^{\pm\pm}(m,n)=\frac{ab}{q^2}\int_0^\infty\int_0^\infty
E(x,y)M^{\pm}_{2i\mu}\left(\frac{4\pi a\sqrt{mx}}{q}\right)
M^{\pm}_{2i\nu}\left(\frac{4\pi b\sqrt{ny}}{q}\right)\,dx\,dy,\]
\[M^+_{2ir}=(4\cosh\pi r)K_{2ir},\quad
M^-_{2ir}=-\frac{\pi}{\cosh\pi r}\{Y_{2ir}+Y_{-2ir}\}.\]
By summing over the residue classes we get
\begin{equation}\label{eq9}
\csillag\sum_{d\smod{q}}\FI_{d/q}= \sum_{\pm\pm}\,\sum_{m,n\geq
1}\lfi(\mp m)\lpsi(\mp n)S_{\overline{\chi\om}}(-h,\pm am\mp
bn;q)E^{\pm\pm}(m,n).
\end{equation}

In order to estimate the twisted Kloosterman sum, we observe that
the greatest common divisor $(-h,\pm am\mp bn,q)$ divides
$N(h,n,a)(h,m,b)$, as follows from the relations $(a,b)=1$ and
$(h,q)=(h,Nab)$. Therefore (\ref{eq10}) and (\ref{eq12}) imply
that
\begin{equation}\label{eq14}
S_{\overline{\chi\om}}(-h,\pm am\mp bn;q)\ll
N^{1/2}(h,m)^{1/2}(h,n)^{1/2}Q^{1/2}(AB)^\ep.
\end{equation}

We estimate $E^{\pm\pm}(m,n)$ by successive applications of
integration by parts and the recurrence relations
\[\frac{d}{dz}\bigl(z^sK_s(z)\bigr)=-z^sK_{s-1}(z),
\quad\frac{d}{dz}\bigl(z^sY_s(z)\bigr)=z^sY_{s-1}(z).\]
Using the first relation we can prove by induction on $k$ that
\[K_s\bigl(\sqrt{z}\bigr)=\sum_{\ka=0}^k c_{\ka k} z^{\ka-\frac{k}{2}}
\left\{K_{s+k}\bigl(\sqrt{z}\bigr)\right\}^{(k)}\]
holds with appropriate constants satisfying
\[c_{\ka k}=c_{\ka k}(s)\ll_k \bigl(1+|s|\bigr)^{k-\ka},\quad 0\leq \ka\leq k.\]
Clearly, for any $\eta>0$ we also have
\[K_s\bigl(\eta\sqrt{z}\bigr)=\eta^{-k}\sum_{\ka=0}^kc_{\ka k}z^{\ka-\frac{k}{2}}
\left\{K_{s+k}\bigl(\eta\sqrt{z}\bigr)\right\}^{(\ka)}.\]
Similarly, for any positive integer $l$ there are
constants
\[d_{\la l}=d_{\la l}(s)\ll_l \bigl(1+|s|\bigr)^{l-\la},\quad 0\leq \la\leq l,\]
such that for any $\theta>0$ we have
\[Y_s\bigl(\theta\sqrt{z}\bigr)=\theta^{-l}\sum_{\la=0}^ld_{\la l}z^{\la-\frac{l}{2}}
\left\{Y_{s+l}\bigl(\theta\sqrt{z}\bigr)\right\}^{(\la)}.\]

By specifying $\eta$ and $\theta$ as
\[\eta=\frac{4\pi a\sqrt{m}}{q},\qquad\qquad\theta=\frac{4\pi b\sqrt{n}}{q},\]
we obtain decompositions of $E^{\pm\pm}(m,n)$ accordingly. In particular, for each pair $(k,l)$ it follows that
\begin{multline}\label{eq103}
E^{\pm\pm}(m,n)\ll_{k,l}
\frac{ab}{q^2}\frac{|\tilde\mu|^{k}|\tilde\nu|^{l}}{\eta^{k}\theta^{l}}
\sup_{M_1,M_2}\sup_{\substack{0\leq\ka\leq k\\0\leq\la\leq l}}\\
\int_0^\infty\int_0^\infty x^{\ka-\frac{k}{2}}y^{\la-\frac{l}{2}}E(x,y)
\left\{M_1\bigl(\eta\sqrt{x}\bigr)\right\}^{(\ka)}
\left\{M_2\bigl(\theta\sqrt{y}\bigr)\right\}^{(\la)}dx\,dy,
\end{multline} 
where
\begin{alignat*}{5}
M_1&\in\Big\{&
(\cosh\pi\mu)&K_{k+2i\mu},\ &
(\cosh\pi\mu)^{-1}&Y_{k+2i\mu},\ &
(\cosh\pi\mu)^{-1}&Y_{k-2i\mu}&
&\Big\},\\
M_2&\in\Big\{&
(\cosh\pi\nu)&K_{l+2i\nu},\ &
(\cosh\pi\nu)^{-1}&Y_{l+2i\nu},\ &
(\cosh\pi\nu)^{-1}&Y_{l-2i\nu}&
&\Big\}.
\end{alignat*}
(\ref{eq109}) shows that each integral above can be rewritten as
\begin{equation}\label{eq104}
\int_{A/a}^{2A/a}\int_{B/b}^{2B/b} 
\left\{x^{\ka-\frac{k}{2}}y^{\la-\frac{l}{2}}E(x,y)\right\}^{(\ka,\la)}
M_1\bigl(\eta\sqrt{x}\bigr)M_2\bigl(\theta\sqrt{y}\bigr)\,dx\,dy.
\end{equation}

We shall pick a pair $(k,l)$ for each $(m,n)$ in such a way, that the following uniform estimates will hold:
\begin{equation}\label{eq105}\begin{split}
M_1\bigl(\eta\sqrt{x}\bigr)&\ll_k|\tilde\mu|^{k+1+\ep}(\eta\sqrt{x}\bigr)^{-1/2},
\qquad x\in[A/a,2A/a];\\
M_2\bigl(\theta\sqrt{y}\bigr)&\ll_l|\tilde\nu|^{l+1+\ep}(\theta\sqrt{y}\bigr)^{-1/2},
\qquad\ \, y\in[B/b,2B/b].
\end{split}\end{equation}
As $|\Re(i\mu)|$ and $|\Re(i\nu)|$ are at most $1/4$, we can refer to the uniform estimates of Section~\ref{bessel} to see that (\ref{eq105}) holds whenever
the assignment $(m,n)\mapsto(k,l)$ is such that
\begin{equation}\label{eq108}\begin{split}
k>0\quad&\Longrightarrow\quad\eta\sqrt{A/a}>1,\\
l>0\quad&\Longrightarrow\quad\theta\sqrt{B/b}>1.
\end{split}\end{equation}
The integral (\ref{eq104}) can be estimated by (\ref{eq7}), (\ref{eq105}), and the relation 
\[\min(A\dd,B\dd)\geq 1,\] 
which follows from (\ref{eq2}). The resulting bound simplifies (\ref{eq103}) to
\[E^{\pm\pm}(m,n)\ll_{k,l}
\frac{AB}{q^2(A+B)}
\frac{|\tilde\mu|^{2k+1+\ep}|\tilde\nu|^{2l+1+\ep}}
{\eta^{k+\frac{1}{2}}\theta^{l+\frac{1}{2}}}
\left(\frac{A}{a}\right)^{-\frac{k}{2}-\frac{1}{4}}
\left(\frac{B}{b}\right)^{-\frac{l}{2}-\frac{1}{4}}
(A\dd)^k(B\dd)^l.\]
Using that $Q\leq q\leq 2Q$ this can be rewritten as
\begin{equation}\label{eq13}
E^{\pm\pm}(m,n)\ll_{k,l} 
\frac{|\tilde\mu\tilde\nu|^\ep(AB)^{1/2}}{\dd Q^2(A+B)}
\left(\frac{A|\tilde\mu|^4(\dd Q)^2}{am}\right)^{\frac{k}{2}+\frac{1}{4}}
\left(\frac{B|\tilde\nu|^4(\dd Q)^2}{bn}\right)^{\frac{l}{2}+\frac{1}{4}}.
\end{equation}
This result is conditional under (\ref{eq105}), but it suggests that
in (\ref{eq9}) we can neglect the contribution of
those pairs $(m,n)$ for which $am/A|\tilde\mu|^4$ or
$bn/B|\tilde\nu|^4$ is greater than
$(\dd Q)^2(AB)^\ep$. 

Indeed, this will be the case if we specify 
$k=\lceil200/\ep\rceil$ or $k=0$ 
(resp. $l=\lceil200/\ep\rceil$ or $l=0$)
depending on whether
$m$ (resp. $n$) 
is large or small in the above sense. We observe that
\[\begin{split}
am>A|\tilde\mu|^4(\dd Q)^2(AB)^\ep\quad&\Longrightarrow\quad\eta\sqrt{A/a}>1,\\
bn>B|\tilde\nu|^4(\dd Q)^2(AB)^\ep\quad&\Longrightarrow\quad\theta\sqrt{B/b}>1,
\end{split}\]
therefore our assignment $(m,n)\mapsto(k,l)$ satisfies (\ref{eq108}),
and with this choice (\ref{eq13}) holds uniformly for all $(m,n)$ with an implied constant depending only on $\ep$.

It follows from (\ref{eq22}) applied to $\phi$ and
$\psi$, that
\begin{equation}\label{eq111}\begin{split}
\sum_{1\leq m\leq x}|\lfi(\mp m)|(h,m)^{1/2}&\ll N^{1/2}|\tilde\mu|^{1/2}x\tau^{1/2}(h),\\
\sum_{1\leq n\leq y}|\lpsi(\mp n)|(h,n)^{1/2}&\ll N^{1/2}|\tilde\nu|^{1/2}y\tau^{1/2}(h).
\end{split}\end{equation}
Combining this bound with (\ref{eq14}) and (\ref{eq13}), we see
that the total contribution to (\ref{eq9}) of small pairs $(m,n)$ is
\[\ll \frac{N^{3/2}|\tilde\mu\tilde\nu|^{3/2+\ep}\dd^3Q^{5/2}}{ab}
\cdot\frac{(AB)^{3/2+\ep}}{A+B}.\]
On the other hand, (\ref{eq14}), (\ref{eq13}) and (\ref{eq111}) similarly
show that the remaining contribution from 
the pairs $(m,n)$ with $m$ or $n$ large is
\[\ll \frac{N^{3/2}|\tilde\mu\tilde\nu|^{3/2+\ep}\dd^3Q^{5/2}}{ab}
\cdot\frac{(AB)^{-50}}{A+B}.\]
To summarize, we have shown that
\[\csillag\sum_{d\smod{q}}\FI_{d/q}
\ll \frac{N^{3/2}|\tilde\mu\tilde\nu|^{3/2+\ep}\dd^3Q^{5/2}}{ab}
\cdot\frac{(AB)^{3/2+\ep}}{A+B}.\]
Hence, by (\ref{eq6}),
\begin{equation}\label{eq16}
\tilde{D}_F(a,b;h)=\frac{1}{2\dd L}\sum_{q\in\CQ}\ \
\csillag\sum_{d\smod{q}}\FI_{d/q}
\ll\frac{N^{3/2}|\tilde\mu\tilde\nu|^{3/2+\ep}\dd^2Q^{3/2}}{ab}
\cdot\frac{(AB)^{3/2+\ep}}{A+B}.
\end{equation}

Inequalities (\ref{eq15}) and (\ref{eq16}) show that the optimal
balance is achieved when
\[\dd^3Q^5\asymp N|\tilde\mu\tilde\nu|(ab)^3.\]
A natural choice is given by
\[\dd^3Q^5=N|\tilde\mu\tilde\nu|(cab)^3,\]
where $c$ is the constant appearing in the remark after
(\ref{eq18}). Then, by (\ref{eq28}),
the conditions of Proposition~\ref{Prop2} are satisfied, that is,
both (\ref{eq3}) and $Q\geq cNab(AB)^\ep$ hold. 
(\ref{eq15}) and (\ref{eq16}) add up to (\ref{eq29}),
using also (\ref{eq30}).

\section{Dyadic decomposition}\label{sect2}

Our aim is to prove Theorem~\ref{Th1} for all test functions
$f(x,y)$ satisfying (\ref{eq1}). We fix an arbitrary smooth
function
\[\rho:(0,\infty)\to\RR\]
whose support lies in $[1,2]$ and which satisfies the following
identity on the positive axis:
\[\sum_{i=-\infty}^\infty\rho\bigl(2^{-i/2}x\bigr)=1.\]
To obtain such a function, we take an arbitrary smooth
$\eta:(0,\infty)\to\RR$ which is constant 0 on $(0,1)$ and
constant 1 on $(\sqrt{2},\infty)$, and then define $\rho$ as
\[\rho(x)=\begin{cases}
\eta(x)&\text{if $0<x\leq\sqrt{2}$,}\\
1-\eta(x/\sqrt{2})&\text{if $\sqrt{2}<x<\infty$.}
\end{cases}\]
According to this partition of unity we decompose $f(x,y)$ as
\[f(x,y)=\sum_{i=-\infty}^\infty\sum_{j=-\infty}^\infty f_{i,j}(x,y),\]
\[f_{i,j}(x,y)=f(x,y)\rho\left(\frac{x}{2^{i/2}X}\right)
\rho\left(\frac{y}{2^{j/2}Y}\right).\] Observe that
\begin{equation}\label{eq8}\text{supp}\,f_{i,j}\subseteq
[A_i,2A_i]\times[B_j,2B_j\bigr],\quad A_i=2^{i/2}X,\quad
B_j=2^{j/2}Y,\end{equation} whence (\ref{eq1}) and $P\geq 1$ show
that
\[\bigl(1+2^{i/2}\bigr)\bigl(1+2^{j/2}\bigr)f_{i,j}^{(k,l)}\ll_{k,l}
A_i^{-k}B_j^{-l}P^{k+l}.\] In other words, the bound (\ref{eq29})
applies uniformly to each function
\[g_{i,j}(x,y)=\bigl(1+2^{i/2}\bigr)\bigl(1+2^{j/2}\bigr)f_{i,j}(x,y)\]
with the corresponding parameters $A=A_i$, $B=B_j$:
\[D_{g_{i,j}}(a,b;h)\ll
P^{11/10}N^{9/5}|\tilde\mu\tilde\nu|^{9/5+\ep}(ab)^{-1/10}(A_i+B_j)^{1/10}(A_iB_j)^{2/5+\ep}.\] This implies, for $\ep<1/10$,
\[D_{f_{i,j}}(a,b;h)\ll
2^{-|i|/5}2^{-|j|/5}P^{11/10}N^{9/5}|\tilde\mu\tilde\nu|^{9/5+\ep}(ab)^{-1/10}(X+Y)^{1/10}(XY)^{2/5+\ep}.\]
Finally,
\[D_f(a,b;h)=\sum_{i=-\infty}^\infty\sum_{j=-\infty}^\infty
D_{f_{i,j}}(a,b;h)\] completes the proof of Theorem~\ref{Th1}.

It should be noted that the trivial upper bound (\ref{eq40})
mentioned in Section~\ref{sect6}
follows by a similar reduction technique from the Cauchy bounds
\[D_{g_{i,j}}(a,b;h)\ll N|\tilde\mu\tilde\nu|^{1/2}(ab)^{-1/2}(A_iB_j)^{1/2}\] 
of Section~\ref{sect3} (cf. (\ref{eq107})).

\section{Bounds for Bessel functions}\label{bessel}

In this section we prove uniform bounds for Bessel functions of the first kind (Proposition~\ref{bessel1}) and of the second and third kinds (Proposition~\ref{bessel2}).

\begin{proposition}\label{bessel1}For any integer $k\geq 1$ the following uniform estimate holds:
\[J_{k-1}(x)\ll\begin{cases}
\frac{x^{k-1}}{2^{k-1}\GG\left(k-\frac{1}{2}\right)},&0<x\leq 1;\\
kx^{-1/2},&1<x.
\end{cases}\]
The implied constant is absolute.
\end{proposition}

\noindent \emph{Proof.} For $x>k^2$ the asymptotic expansion of $J_{k-1}$ 
(see Section~7.21 of \cite{Wat})
provides the stronger estimate $J_{k-1}(x)\ll x^{-1/2}$
with an absolute implied constant.

For $1<x\leq k^2$ we use Bessel's original 
integral representation (see Section~2.2 of \cite{Wat}),
\[J_{k-1}(x)=\frac{1}{2\pi}\int_0^{2\pi}\cos\bigl((k-1)\theta-x\sin\theta\bigr)\,d\theta,\]
to deduce that in this range 
\[|J_{k-1}(x)|\leq 1\leq kx^{-1/2}.\]

For the remaining range $0<x\leq 1$ the required estimate follows
from the Poisson-Lommel integral representation (see Section~3.3 of \cite{Wat})
\[J_{k-1}(x)=\frac{x^{k-1}}{2^{k-1}\GG\left(k-\frac{1}{2}\right)\GG\left(\frac{1}{2}\right)}
\int_0^\pi\cos(x\cos\theta)\sin^{2k-2}\theta\,d\theta.\qed\]

\begin{proposition}\label{bessel2}For any $\si>0$ and $\ep>0$ 
the following uniform estimates hold in the strip
$|\Re s|\leq\si$:
\[\begin{split}e^{-\pi|\Im s|/2}Y_s(x)
&\ll\begin{cases}
\bigl(1+|\Im s|\bigr)^{\si+\ep}x^{-\si-\ep},&0<x\leq 1+|\Im s|;\\
\bigl(1+|\Im s|\bigr)^{-\ep}x^{\ep},&1+|\Im s|<x\leq 1+|s|^2;\\
x^{-1/2},&1+|s|^2<x.
\end{cases}\\\\
e^{\pi|\Im s|/2}K_s(x)
&\ll\begin{cases}
\bigl(1+|\Im s|\bigr)^{\si+\ep}x^{-\si-\ep},&0<x\leq 1+\pi|\Im s|/2;\\
e^{-x+\pi|\Im s|/2}x^{-1/2},&1+\pi|\Im s|/2<x.\\
\end{cases}\end{split}\]
The implied constants depend only on $\si$ and $\ep$.
\end{proposition}

\noindent \emph{Proof.} The last estimate for $Y_s$ follows from its asymptotic expansion (see Section~7.21 of \cite{Wat}). The last estimate for $K_s$ follows from Schl\"afli's integral
representation (see Section~6.22 of \cite{Wat}),
\[K_s(x)=\int_0^\infty e^{-x\cosh t}\cosh st.\,dt,\]
by noting that 
\[\cosh t\geq 1+t^2/2\quad\text{and}\quad|\cosh st|\leq e^{\si t}.\]

We shall deduce the remaining uniform bounds from the
integral representations
\[\begin{split}
4K_s(x)(x)&=\frac{1}{2\pi i}\int_{\cC}
\GG\left(\frac{w-s}{2}\right)\GG\left(\frac{w+s}{2}\right)
\left(\frac{x}{2}\right)^{-w}\,dw,\\\\
-2\pi Y_s(x)(x)&=\frac{1}{2\pi i}\int_{\cC}
\GG\left(\frac{w-s}{2}\right)\GG\left(\frac{w+s}{2}\right)
\cos\left(\frac{\pi}{2}(w-s)\right)
\left(\frac{x}{2}\right)^{-w}\,dw,
\end{split}\]
where the contour $\cC$ is a broken line of 2 infinite and 3 finite segments joining the points \[\begin{matrix}
-\ep-i\infty,
&-\ep-i\bigl(2+|\Im s|\bigr),
&\si+\ep-i\bigl(2+|\Im s|\bigr),\\
\si+\ep+i\bigl(2+|\Im s|\bigr),
&-\ep+i\bigl(2+|\Im s|\bigr),
&-\ep+i\infty.
\end{matrix}\]
These formulae follow by analytic continuation from the well-known but more restrictive inverse Mellin transform representations of the $K$- and $Y$-Bessel functions, cf. formulae 6.8.17 and 6.8.26 in \cite{Er}.

If we write in the second formula
\[\cos\left(\frac{\pi}{2}(w-s)\right)=
\cos\left(\frac{\pi}{2}w\right)\cos\left(\frac{\pi}{2}s\right)+
\sin\left(\frac{\pi}{2}w\right)\sin\left(\frac{\pi}{2}s\right),\]
then it becomes apparent that the remaining inequalities of the lemma can be deduced from the uniform bound
\begin{multline*}\int_{\cC}e^{\pi\max(|\Im s|,|\Im w|)/2}
\left|\GG\left(\frac{w-s}{2}\right)\GG\left(\frac{w+s}{2}\right)
\left(\frac{x}{2}\right)^{-w}\,dw\right|\\
\ll_{\si,\ep}
\left(\frac{x}{1+|\Im s|}\right)^{-\si-\ep}+\left(\frac{x}{1+|\Im s|}\right)^{\ep}.
\end{multline*}
By introducing the notation
\[G(s)=e^{-\pi|\Im s|/2}\GG(s),\]
\[M_s(x)=\int_{\cC}
\left|G\left(\frac{w-s}{2}\right)G\left(\frac{w+s}{2}\right)
\left(\frac{x}{2}\right)^{-w}\,dw\right|,\]
the previous inequality can be rewritten as
\begin{equation}\label{eq112}
M_s(x)\ll_{\si,\ep}
\left(\frac{x}{1+|\Im s|}\right)^{-\si-\ep}+\left(\frac{x}{1+|\Im s|}\right)^{\ep}.
\end{equation}

\begin{Case} $|\Im s|\leq 1$.

If $w$ lies on either horizontal segments of $\cC$ or on the finite vertical segment joining $\si+\ep\pm i\bigl(2+|\Im s|\bigr)$, then $w\pm s$ varies in a fixed compact set (depending only on $\si$ and $\ep$) disjoint from the negative axis $\Re z\leq 0$. It follows that for these values $w$ we have
\[G\left(\frac{w-s}{2}\right)G\left(\frac{w+s}{2}\right)\ll_{\si,\ep} 1,\]
i.e.,
\[G\left(\frac{w-s}{2}\right)G\left(\frac{w+s}{2}\right)\left(\frac{x}{2}\right)^{-w}
\ll_{\si,\ep}x^{-\si-\ep},\]
and the same bound holds for the contribution of these values to $M_s(x)$.

If $w$ lies on either infinite vertical segments of $\cC$, then
\[|\Im(w\pm s)|\asymp|\Im w|>1,\]
whence Stirling's approximation yields
\[G\left(\frac{w-s}{2}\right)G\left(\frac{w+s}{2}\right)\asymp_\ep|\Im w|^{-\ep-1}.\]
It follows that the contribution of the infinite segments to $M_s(x)$ is 
$\ll_{\si,\ep}x^\ep$.

Altogether we infer that
\[M_s(x)\ll_{\si,\ep}x^{-\si-\ep}+x^\ep,\]
which is equivalent to (\ref{eq112}).
\end{Case}

\begin{Case} $|\Im s|>1$.

If $w$ lies on either horizontal segments of $\cC$, then
\[|\Im(w\pm s)|\asymp|\Im s|,\]
whence Stirling's approximation yields
\[G\left(\frac{w-s}{2}\right)G\left(\frac{w+s}{2}\right)\asymp_{\si,\ep}|\Im s|^{\Re w-1},\]
i.e.,
\[G\left(\frac{w-s}{2}\right)G\left(\frac{w+s}{2}\right)\left(\frac{x}{2}\right)^{-w}\asymp_{\si,\ep}\frac{1}{|\Im s|}\left(\frac{|\Im s|}{x}\right)^{\Re w}.\]
It follows that the contribution of the horizontal segments to $M_s(x)$ is
\[\ll_{\si,\ep}
|\Im s|^{-1+\si+\ep}x^{-\si-\ep}+|\Im s|^{-1-\ep}x^{\ep}.\]

If $w$ lies on the finite vertical segment of $\cC$ 
joining $\si+\ep\pm i\bigl(2+|\Im s|\bigr)$, then
\[\Re(w\pm s)\geq\ep
\qquad\text{and}\qquad
\max|\Im(w\pm s)|\asymp|\Im s|,\]
whence Stirling's approximation implies
\[G\left(\frac{w-s}{2}\right)G\left(\frac{w+s}{2}\right)\ll_{\si,\ep}
\begin{cases}
|\Im s|^{\si+\ep/2-1/2}&\text{if \ $\min|\Im(w\pm s)|\leq 1$;}\\
|\Im s|^{\si+\ep-1}&\text{if \ $\min|\Im(w\pm s)|>1$.}
\end{cases}\]
It follows that the contribution of the finite vertical segment to $M_s(x)$ is
\[\ll_{\si,\ep}|\Im s|^{\si+\ep}x^{-\si-\ep}.\]

If $w$ lies on either infinite vertical segments of $\cC$, then
\[|\Im(w\pm s)|\asymp|\Im w|>|\Im s|,\]
whence Stirling's approximation yields
\[G\left(\frac{w-s}{2}\right)G\left(\frac{w+s}{2}\right)\asymp_\ep|\Im w|^{-\ep-1}.\]
It follows that the contribution of the infinite vertical segments to $M_s(x)$ is
\[\ll_{\si,\ep}|\Im s|^{-\ep}x^{\ep}.\]

Altogether we infer that
\[M_s(x)\ll_{\si,\ep}|\Im s|^{\si+\ep}x^{-\si-\ep}+|\Im s|^{-\ep}x^\ep,\]
which is equivalent to (\ref{eq112}).

\end{Case}

The proof of Proposition~\ref{bessel2} is complete.\qed
\chapter{Twists of Maass forms: a subconvex bound for $L$-functions}\label{subconvex}

\section{Overview}

We shall prove a subconvex estimate on the critical
line for $L$-functions associated to character twists of a fixed
holomorphic or Maass cusp form $\phi$ of arbitrary level and
nebentypus. We borrow notation from Section~\ref{sect1}, and we also
refer the reader to Section~\ref{size} for an introduction. The result, in less explicit form, will also appear in \cite{Ha2}.

We assume that $\phi$ is a \emph{primitive form}, that is, 
a newform in the sense of \cite{AL,Li,ALi} normalized so that
${\rho}_\phi(1)=1$. If we renormalize the Fourier coefficients of $\phi$
as
\[\lfi(n)=|n|^{\frac{1-k}{2}}{\rho}_\phi(n),\]
then $\lfi(n)$ $(n\geq 1)$ defines a
character of the corresponding Hecke algebra, while
$\lfi(-n)=\pm\lfi(n)$ (with a constant sign) when $\phi$ is a
Maass form. In other words, $\phi$ defines a cuspidal automorphic
representation of $\GL_2$ over $\QQ$ with arithmetic conductor
$N$. The contragradient
representation corresponds to the primitive cusp form
$\tilde\phi(z)=\bar\phi(-\bar z)$ with renormalized Fourier coefficients
$\lambda_{\tilde\phi}(n)=\blfi(n)$. We note that by the powerful results
of Iwaniec \cite{I2} and Hoffstein--Lockhart \cite{HL}, the old normalization
(\ref{eq200}) and the present one are essentially the same in that the scaling factor $c$ between them satisfies
\[N^{-\ep}|\tilde\mu|^{-\ep}\ll_\ep c
\ll_\ep N^\ep|\tilde\mu|^{\ep}.\]

We consider the twisted
representations $\phi\otimes\chi$ as $\chi$ runs through the automorphic
representations of $\GL_1$ over $\QQ$, that is, the primitive Dirichlet
characters of the rational integers. 
In order to simplify our discussion, we shall
assume that $q$, the conductor of $\chi$, 
is prime to $N$. Then the analytic conductor of
$\phi\otimes\chi$ satisfies
\begin{equation}\label{eq202}
C(s,\phi\otimes\chi)\asymp q^2 N\bigl(|s|^2+\tmu^2\bigr),
\qquad\Re s=\frac{1}{2},\end{equation}
and for $\Re s>1$ the associated $L$-function is given by
\[L(s,\phi\otimes\chi)=\sum_{n=1}^\infty\frac{\lfi(n)\chi(n)}{n^s}.\]
For a fixed point $s$ on the critical line the convexity bound (\ref{eq80})
implies that
\[L(s,\phi\otimes\chi)\ll_{\ep}
|s|^{1/2+\ep}N^{1/4+\ep}\tmu^{1/2+\ep}q^{1/2+\ep}.\]
Our aim is to decrease the exponent $1/2$ of $q$ 
and still maintain polynomial
control in the other parameters $|s|$, $N$, $\tmu$.

\begin{theorem}\label{Th2} Suppose that $\phi$ is a
primitive holomorphic or Maass cusp form of Archimedean size $\tmu$, 
level $N$ and arbitrary
nebentypus character mod $N$. Let $\Re s=1/2$ and 
$q$ be an integer prime to $N$. 
If $\chi$ is a primitive Dirichlet character modulo $q$,
then
\begin{equation}\label{eq19}
L(s,\phi\otimes\chi)\ll |s|^{1+\ep}N^{9/8+\ep}\tmu^{27/20+\ep}
q^{1/2-1/54+\ep},
\end{equation}
where the implied constant depends only on $\ep$.
\end{theorem}

A similar estimate with $q$-exponent $1/2-1/22$ was proved for holomorphic
forms of full level in \cite{DFI2}, and the improved exponent
$1/2-7/130$ follows for holomorphic forms of arbitrary level as
a special case of the main result in \cite{Sa3}.
Duke, Friedlander and Iwaniec anticipated their
method to be extendible to more general $L$-functions of rank two,
and the present chapter is indeed an extension of their work.
The very general Vorono\"\i-formula of Michel enables one to establish 
Theorem~\ref{Th2} in slightly stronger form, e.g. with the original 
$q$-exponent $1/2-1/22$ of \cite{DFI2}. See \cite{Mi2} for details.

Combining the estimate (\ref{eq19}) at the central point $s=1/2$
with Waldspurger's theorem \cite{Wa} (see also \cite{K,Sh}), we get
the bound
\[c(q)\ll_\ep q^{1/4-1/108+\ep},\quad\text{$q$ square-free}\]
for the normalized Fourier coefficients of half-integral weight
forms of arbitrary level. Such a nontrivial bound is the key step
in the solution of the general ternary Linnik problem given by
Duke and Schulze-Pillot \cite{D,D-SP}.

The proof of Theorem~\ref{Th2} is presented in
Sections~\ref{sect4} through \ref{sect5}. In Section~\ref{sect4} we reduce
(\ref{eq19}), via the approximate functional equation of 
Chapter~\ref{functional}, to an
inequality about certain finite sums involving at most
$C(s,\phi\otimes\chi)^{1/2+\ep}$ terms (cf. (\ref{eq202})). 
We prove this inequality in Section~\ref{sect20} by employing
the amplification method.
As discussed in Section~\ref{amplification}, the idea is to consider a
suitably weighted second moment of the finite sums arising from the family
$\phi\otimes\chi$ of cusp forms ($\chi$ varies, $\phi$ is fixed).
We choose the weights (called amplifiers) in such a way that one of the characters
$\chi$ is emphasized, while the second moment average is still of moderate size.
This forces, by positivity, $L(s,\phi\otimes\chi)$ to be small. In
the course of evaluating the amplified second moment we
encounter diagonal and off-diagonal terms. The off-diagonal terms 
decompose to shifted convolution sums, and at this point we apply
Theorem~\ref{Th1}.

\section{Approximate functional equation}\label{sect4}

Using the approximate functional equation in the form Corollary~\ref{Cor1},
we can see that (\ref{eq19}) is equivalent to
\[\sum_{n\leq C^{1/2+\ep}}\frac{\lfi(n)\chi(n)}{n^s}f\left(\frac{n}{\sqrt{C}}\right)
\ll_\ep |s|^{1+\ep}N^{9/8+\ep}\tmu^{27/20+\ep}q^{1/2-1/54+\ep},\]
where
\[C=C(s,\phi\otimes\chi)\ll |s|^2N\tmu^2q^2,\]
and $f:(0,\infty)\to\CC$ is a smooth function satisfying (\ref{eq05}) and (\ref{eq012}) with $m=2$. In particular, we can write the left hand side as
\[\sum_{n\leq C^{1/2+\ep}}\frac{\lfi(n)\chi(n)g(n)}{\sqrt{n}},\]
where
\[g(x)=x^{1/2-s}f\left(\frac{x}{\sqrt{C}}\right)\]
satisfies the uniform bounds
\[g^{(k)}(x)\ll_k |s|^kx^{-k}.\]

Therefore, applying partial summation and a smooth dyadic
decomposition, we can reduce Theorem~\ref{Th2} to the following
\begin{proposition}\label{Prop4}
Let $1\leq T\leq\bigl(|s|N^{1/2}\tmu q\bigr)^{1+\ep}$ 
and $W$ be a smooth complex valued function supported
in $[T,2T]$ such that $W^{(k)}\ll_k |s|^k T^{-k}$. Then
\[\sum_{n=1}^\infty\lfi(n)\chi(n)W(n)\ll 
|s|^{5/6+\ep}N^{25/24+\ep}\tmu^{71/60+\ep}q^{17/54+\ep}T^{2/3},\]
where the implied constant depends only on $\ep$.
\end{proposition}

\section{Amplification}\label{sect20}

Our purpose is to prove Proposition~\ref{Prop4}. As in \cite{DFI2},
we shall estimate from both ways the amplified second moment
\[S=\csillag\sum_{\om\bmod{q}}
\left|\sum_{1\leq l\leq L}\bar\chi(l)\om(l)\right|^2|S_\om|^2,\]
where $\om$ runs through the primitive characters modulo $q$, $L$
is a parameter to be chosen later in terms of $M$ and $q$, and
\[S_\om=\sum_{n=1}^\infty\lfi(n)\om(n)W(n).\]
Assuming $L\geq c(\ep)q^\ep$, it follows, using the
result of Jacobsthal \cite{Ja} that the largest gap between
reduced residue classes mod $q$ is of size $\ll q^\ep$, that
\begin{equation}\label{eq20}
S\gg q^{-\ep}L^2|S_\chi|^2.
\end{equation}
Here and in the sequel implied contants may depend on $\ep$.

On the other hand, expanding each primitive $\om$ in $S$ using
Gauss sums and then extending the resulting summation to all
characters mod $q$, we get by orthogonality,
\[S\leq\frac{\vfi(q)}{q}\csillag\sum_{d\smod{q}}
\left|\sum_m a(m)e_q(dm)\right|^2,\] where
\[a(m)=\sum_{\substack{ln=m\\1\leq l\leq
L}}\bar\chi(l)\lfi(n)W(n).\] It is clear that the coefficients
$a(m)$ are supported in the interval $[1,M]$, where $M=2LT$.
Extending the summation to all residue classes $d$, the previous
inequality becomes
\begin{equation}\label{eq21}
S\leq\vfi(q)\sum_{h\equiv 0\smod{q}}D(h),
\end{equation}
where
\[D(h)=\sum_{m_1-m_2=h}a(m_1)\bar a(m_2).\]

We estimate the diagonal contribution $D(0)$ using the 
following Rankin--Selberg bound (Theorem~8.3 in \cite{I}):
\[\sum_{1\leq n\leq x}|\lfi(n)|^2\ll N^\ep|\tilde\mu|^\ep x.\]
Indeed, by $W\ll 1$ we
get
\[\begin{split}
D(0)&=\sum_m|a(m)|^2\ll\sum_{\substack{l_1n_1=l_2n_2\\1\leq
l_1,l_2\leq L\\T\leq n_1,n_2\leq 2T}}\lfi(n_1)\blfi(n_2)\\\\
&\ll\sum_{\substack{1\leq l\leq L\\T\leq n\leq
2T}}|\lfi(n)|^2\tau(nl)\ll N^\ep\tmu^\ep M^\ep L\sum_{T\leq n\leq
2T}|\lfi(n)|^2,
\end{split}\]
whence
\begin{equation}\label{eq23}
D(0)=\sum_m|a(m)|^2\ll N^\ep\tmu^\ep M^{1+\ep}.
\end{equation}

We estimate the non-diagonal terms $D(h)$ $(h\neq 0)$ using
Theorem~\ref{Th1}. Clearly, we can rewrite each
term as
\[D(h)=\sum_{1\leq l_1,l_2\leq L}\bar\chi(l_1)\chi(l_2)
\sum_{l_1n_1-l_2n_2=h}\lfi(n_1)\blfi(n_2)W(n_1)\bar W(n_2).\] The
inner sum is of type (\ref{eq24}), because $\blfi(n)$ is just the
$n$-th renormalized Fourier coefficient of the contragradient cusp
form $\tilde\phi(z)=\bar\phi(-\bar z)$. For each pair $(l_1,l_2)$
we apply Theorem~\ref{Th1} with $a=l_1/(l_1,l_2)$,
$b=l_2/(l_1,l_2)$, $P=2|s|$, $X=aT$ and $Y=bT$ to conclude that
\begin{equation}\label{eq25}\begin{split}
D(h)
&\ll L^2|s|^{11/10}N^{9/5+\ep}\tmu^{9/5+\ep}
(a+b)^{1/10}(ab)^{3/10+\ep}T^{9/10+\ep}\\
&\ll |s|^{11/10}N^{9/5+\ep}\tmu^{9/5+\ep}L^{27/10+\ep}T^{9/10+\ep}.
\end{split}\end{equation}

\section{Optimizing parameters}\label{sect5}

Inserting the bounds (\ref{eq23}) and (\ref{eq25}) into
(\ref{eq21}), it follows that
\[S\ll N^\ep\tmu^\ep M^\ep\vfi(q)\left(M+\frac{M}{q}
|s|^{11/10}N^{9/5}\tmu^{9/5}L^{27/10}T^{9/10}\right).\]
This shows that the optimal choice for $L$ is provided by
\[q\asymp L^{27/10}T^{9/10}.\]
In order to maintain $L\geq c(\ep)q^\ep$, we choose
\begin{equation}\label{eq26}
\bigl(|s|N^{1/2}\tmu\bigr)^{9/10+\ep}q=L^{27/10}T^{9/10}.
\end{equation}
This shows that
\[S\ll |s|^{2+\ep}N^{9/4+\ep}\tmu^{27/10+\ep}qM^{1+\ep},\]
and then (\ref{eq20}) yields
\[S_\chi\ll q^\ep L^{-1}|S|^{1/2}\ll 
\bigl(|s|^{2}N^{9/4}\tmu^{27/10}qT/L)^{1/2+\ep}.\] Substituting
(\ref{eq26}) we get
\[S_\chi\ll 
\left\{|s|^{2}N^{9/4}\tmu^{27/10}qT
(|s|N^{1/2}\tmu\bigr)^{-1/3}
q^{-10/27}T^{1/3}\right\}^{1/2+\ep},\]
which is precisely the conclusion of Proposition~\ref{Prop4}.\qed 

The proof of Theorem~\ref{Th2} is complete.

\chapter{Shifted convolution sums and spectral theory}\label{convolution}

\section{Overview}

We shall obtain a fairly precise description of the continuous span of the
functions $H_{s,0,i\mu}$ corresponding to values $s$ on a vertical line 
$\si+i\RR$, $\si>1$. These functions play an important role in the Sarnak--Selberg
spectral method applied to Maass forms. We refer the reader to Section~\ref{shifted2} for an introduction. For convenience we shall assume that $\mu\in\RR$.

By definition,
\[\begin{split}H_{s,0,i\mu}(u)
&=\int_0^\infty\tilde W_{0,i\mu}\bigl(|u+1|y\bigr)\bar{\tilde W}_{0,i\mu}\bigl(|u-1|y\bigr)y^{s-2}\,dy\\\\
&=\frac{|u^2-1|^\frac{1}{2}}{\pi}
\int_0^\infty K_{i\mu}\left(\frac{|u+1|y}{2}\right)
K_{i\mu}\left(\frac{|u-1|y}{2}\right)y^{s}\frac{dy}{y}.
\end{split}\]
In particular, $H_{s,0,i\mu}(u)$ is an even function of $u$, therefore we can regard it as a function on the positive axis.
Combining formulae 6.576.4 and 9.134.3 from \cite{GR}, we can see that
\[H_{s,0,i\mu}(u)=M(s)u|1-u^{-2}|^{\frac{1}{2}+i\mu}G_s(u),\]
where
\[M(s)=\frac{2^{2s-3}\GG\left(\frac{s}{2}-i\mu\right)
\GG^2\left(\frac{s}{2}\right)\GG\left(\frac{s}{2}+i\mu\right)}{\pi\GG(s)},\]
and
\[G_s(u)=\begin{cases}
u^{2i\mu}F\left(\tfrac{s}{2}+i\mu,\tfrac{1}{2}+i\mu;\tfrac{s}{2}+\tfrac{1}{2};u^2\right)
,&\quad\text{$0\leq u<1$;}\\
u^{-s}F\left(\tfrac{s}{2}+i\mu,\tfrac{1}{2}+i\mu;\tfrac{s}{2}+\tfrac{1}{2};u^{-2}\right)
,&\quad\text{$1<u$.}
\end{cases}\]

This explicit decomposition reduces our task to analyze the set of functions $V$
on the positive axis that can be represented in the form
\begin{equation}\label{eq95}
V(u)=\frac{1}{2\pi i}\int_{(\si)}V^\jb(s)G_s(u)\,ds.\end{equation}
In our formal definition we include an assumption on the growth rate
of $V^\jb(s)$, which is natural and justified by the dependencies on $s$
of the bounds in Lemmata \ref{F} and \ref{F2}.

\begin{Def} Let $V$ be an arbitrary complex valued function on 
the positive axis $(0,\infty)$, and $V^\jb(s)$ be a complex valued
function on the vertical line $\si+i\RR$ ($\si>1$), such that
\begin{equation}\label{eq221}
\int_{(\si)}|s|^{3/2+\ep}|V^\jb(s)|\,ds<\infty\end{equation}
holds for some $\ep>0$. 
Then $V^\jb$
is a $\jb$ transform of $V$ if (\ref{eq95}) is valid for all
$u>0$, $u\neq 1$.\end{Def}

\begin{theorem}\label{Th91}Suppose that an arbitrary function
$V:(0,\infty)\to\CC$ has a $\jb$ transform on the vertical line $\si+i\RR$ ($\si>1$). Then $V$ is continuous at all points $u\neq 1$,
the Mellin transform $V^*(z)$ of $V$ is defined in
$0<\Re z<1$, and $\frac{\GG\left(\frac{z}{2}+\frac{1}{2}\right)}
{\GG\left(\frac{z}{2}+i\mu\right)}V^*(z)$ extends to a bounded holomorphic
function in every half-plane $\Re z<\si_0<\si$. Conversely, let $V:(0,\infty)\to\CC$
be an arbitrary function which is continuous at all points $u\neq 1$ and
has Mellin transform $V^*(z)$ defined in 
$0<\Re z<1$. If $K(z)=\frac{\GG\left(\frac{z}{2}+\frac{1}{2}\right)}
{\GG\left(\frac{z}{2}+i\mu\right)}V^*(z)$ extends to a holomorphic
function in some half-plane $\Re z<\si_0$ ($\si_0>1$) satisfying
$K(z)\ll(1+|z|)^{-A}$ for some $A>2$, then $V$ has a $\jb$ transform
$V^\jb(s)$, which extends to a holomorphic function
in $0<\Re s<\si_0$ satisfying $V^\jb(s)\ll_{\si,A}\bigl(1+|s|\bigr)^{-A-1/2}$.\end{theorem}

The theorem shows that the functions $H_{s,0,i\mu}$ form an incomplete 
system in the sense that some of the very natural functions $V$
are excluded from their continuous span. For the Sarnak--Selberg method this negative
conclusion has the message that
the Maass operators must play a crucial role in a successful analysis.

\begin{Cor}\label{Cor91}Let $V:(0,\infty)\to\CC$ be an arbitrary
function compactly supported 
in $(0,1)\cup(1,\infty)$. If $V$ has a $\jb$ transform,
then it is identically zero.\end{Cor}

\begin{proof}By Theorem~\ref{Th91}, $V$ is a continuous function
of compact support whose Mellin transform vanishes at all negative 
odd integers. In other words, $V(u)$ is orthogonal to all functions
$u^{-2k}$ $(k=1,2,\dots)$. It follows that $u^{-2}V(u)$ is 
orthogonal to all functions $p(u^{-2})$, where $p$ is an arbitrary
complex polynomial. These functions are dense among continuous
functions on a compact interval
by Weierstrass' approximation theorem, hence $V=0$.\end{proof}

\begin{Cor}\label{Cor92}Let $V:(0,\infty)\to\CC$ be an arbitrary function.
If $V(u/c)$ has a $\jb$ transform for every $c>0$, then $V$
is identically zero.\end{Cor}

\begin{proof}By Theorem~\ref{Th91}, 
\[K(z)=\frac{\GG\left(\frac{z}{2}+\frac{1}{2}\right)}
{\GG\left(\frac{z}{2}+i\mu\right)}V^*(z)\]
is defined in the strip $0<\Re z<1$ and extends to
a holomorphic function in the half-plane $\Re z<1$.
Moreover, for any $c>0$, $c^zK(z)$ is bounded. This forces $K(z)=0$ 
as follows. Let $f(u)$ be the inverse Mellin
transform of $K(z)/(2-z)^2$, i.e.,
\[f(u)=\frac{1}{2\pi i}\int_{(\si)}u^{-z}\frac{K(z)}{(2-z)^2}\,dz\]
for any $\si<1$. We can see that $f(u)$ is independendent of
the particular line of integration. However, the assumption that
$c^zK(z)$ is bounded for any $c>0$ implies uniform bounds of the form
\[f(u)\ll_c(cu)^{-\si},\quad u>0,\quad \si<1,\]
the implied constant depending on $c$ only. By letting $\si\to -\infty$, 
we can conclude, for each $c>0$, that $f(u)$ vanishes on $(1/c,\infty)$.
Hence $f(u)$ is identically zero, and
\[K(z)=(2-z)^2\int_0^\infty u^zf(u)\frac{du}{u}=0,\quad\Re z<1,\]
as claimed. 

Therefore $V^*$ must vanish in $0<\Re z<1$,
which shows that $V(u)=0$ as long as $u\neq 1$. We can repeat the argument
with $V(2u)$ in place of $V(u)$ to see that $V(1)=0$ must hold
as well.\end{proof}

\section{The integral transform}

In this section we prove Theorem~\ref{Th91}.
To prove the first part, we shall assume
that (\ref{eq95}) holds for all $u>0$, $u\neq 1$, where
$\si>1$, and $V^\jb$ is a complex valued function on the vertical line
$\si+i\RR$ satisfying (\ref{eq221}) 
for some $\ep>0$.
By applying formally the Mellin transform on both sides, we get
\begin{equation}\label{eq96}
V^*(z)=\frac{1}{2\pi
i}\int_{(\si)}V^\jb(s)G_s^*(z)\,ds.\end{equation} This step is
justified by Fubini's theorem, if a sufficient uniform bound is
provided for $G_s(u)$. We need to give a uniform estimate for the
hypergeometric functions appearing in $G_s(u)$.

\begin{Lemma}\label{F}Let $\si>1$ and $0\leq u<1$. Then for any $\ep>0$
the following uniform bound holds on the vertical line $\si+i\RR$:
\begin{equation}\label{eq97}
F\left(\tfrac{s}{2}+i\mu,\tfrac{1}{2}+i\mu;\tfrac{s}{2}+\tfrac{1}{2};u\right)
\ll|s|^{1/2+\ep}.\end{equation}
The implied constant depends only on $\si$ and $\ep$.
\end{Lemma}

We postpone the proof of this lemma to the next section. It shows
that the Mellin transforms $V^*(z)$ and  $G_s^*(z)$ exist
for $0<\Re z<\si$, and that (\ref{eq96}) is valid in this strip.

We can also see that $V$ is continuous at $u\neq 1$, because the functions
$G_s(u)$ are sufficiently uniformly continuous at these points. The relevant
estimate reads as follows.

\begin{Lemma}\label{F2}Let $\si>1$ and $0\leq v<u<1$. Then for any $\ep>0$ 
the following uniform bound holds on the vertical line $\si+i\RR$:
\[F\left(\tfrac{s}{2}+i\mu,\tfrac{1}{2}+i\mu;\tfrac{s}{2}+\tfrac{1}{2};u\right)-
F\left(\tfrac{s}{2}+i\mu,\tfrac{1}{2}+i\mu;\tfrac{s}{2}+\tfrac{1}{2};v\right)
\ll(u-v)|s|^{3/2+\ep}.\]
\end{Lemma}

We omit the proof of this result, as it is almost identical to that of Lemma~\ref{F}.

It is essential that $G_s^*(z)$ can be determined explicitly.
It is given by a special case of formula 2.21.1.3 from \cite{Pr}:
\[\begin{split}
\tfrac{\GG(\al)\GG(a-\al)}{\GG(1-b+\al)\GG(c-\al)}
&=\tfrac{\GG(a)}{\GG(1-b)\GG(c)}\int_0^1u^\al
F(a,b;c;u)\tfrac{\,du}{u}\\\\
&+\tfrac{\GG(a)}{\GG(c-a)\GG(a-b+1)}\int_1^\infty u^{\al-a}
F\left(a,a-c+1;a-b+1;\tfrac{1}{u}\right)\tfrac{\,du}{u}.\end{split}\]

\medskip\noindent The formula is valid as long as $\Re(c-a-b)>-1$,
$0<\Re\al<\Re a$ and all the gamma values are finite on the right
hand side. It can be verified formally by regarding the left hand
side as a function of $\al$ and evaluating its inverse Mellin
transform in $u$. By specializing the above formula to
\[\al=\tfrac{z}{2}+i\mu,\qquad a=\tfrac{s}{2}+i\mu
,\qquad b=\tfrac{1}{2}+i\mu,\qquad c=\tfrac{s}{2}+\tfrac{1}{2},\]
and replacing $u$ by $u^2$ in both integrals, we obtain the
following result.

\begin{Lemma}\label{G}For $0<\Re z<\si$ the Mellin transform of
$G_s(u)$ is given by
\[G_s^*(z)=c_{i\mu}
\frac{\GG\left(\frac{s}{2}+\frac{1}{2}\right)}{\GG\left(\frac{s}{2}+i\mu\right)}
\frac{\GG\left(\frac{s-z}{2}\right)}{\GG\left(\frac{s-z}{2}+\frac{1}{2}-i\mu\right)}
\frac{\GG\left(\frac{z}{2}+i\mu\right)}{\GG\left(\frac{z}{2}+\frac{1}{2}\right)},\]
where $c_{i\mu}$ abbreviates the constant
$\tfrac{1}{2}\GG\left(\tfrac{1}{2}-i\mu\right)$.
\end{Lemma}

In particular, for $0<\Re z<\si$, (\ref{eq96}) can be
rewritten as
\begin{equation}\label{eq98}
\frac{\GG\left(\frac{z}{2}+\frac{1}{2}\right)}{\GG\left(\frac{z}{2}+i\mu\right)}
V^*(z)=\frac{1}{2\pi i}\int_{(\si)}V^\jb(s)
\frac{\GG\left(\frac{s}{2}+\frac{1}{2}\right)}{\GG\left(\frac{s}{2}+i\mu\right)}
\frac{c_{i\mu}\GG\left(\frac{s-z}{2}\right)}{\GG\left(\frac{s-z}{2}+\frac{1}{2}-i\mu\right)}
\,ds.\end{equation} The 
right hand side of this equation defines a bounded holomorphic function in
every half-plane $\Re z<\si_0<\si$, which concludes the proof of the first part of the theorem.

We turn to the second part of Theorem~\ref{Th91}. Let $V:(0,\infty)\to\CC$
be an arbitrary function which is continuous at all points $u\neq 1$.
We assume that the Mellin transform $V^*(z)$ is defined for
$0<\Re z<1$ such that $\frac{\GG\left(\frac{z}{2}+\frac{1}{2}\right)}
{\GG\left(\frac{z}{2}+i\mu\right)}V^*(z)$ extends to a holomorphic
function in some half-plane $\Re z<\si_0$ ($\si_0>1$). If, in addition, we have a uniform bound
\begin{equation}\label{eq220}
\frac{\GG\left(\frac{z}{2}+\frac{1}{2}\right)}
{\GG\left(\frac{z}{2}+i\mu\right)}V^*(z)\ll(1+|z|)^{-A},\quad\Re z<\si_0
\end{equation}
for some $A>2$, then the inverse Mellin transform of the left hand side
is a continuous function $k:(0,\infty)\to\CC$
vanishing on $(0,1)$ such that
\begin{equation}\label{eq912}k^*(z)=\frac{\GG\left(\frac{z}{2}+\frac{1}{2}\right)}{\GG\left(\frac{z}{2}+i\mu\right)}
V^*(z).\end{equation} 

We claim that there is a continuous function
$l:(0,\infty)\to\CC$ vanishing on $(0,1)$ such that (\ref{eq98}) is solved by
\begin{equation}\label{eq911}
V^\jb(s)=\frac{\GG\left(\frac{s}{2}+i\mu\right)}
{\GG\left(\frac{s}{2}+\frac{1}{2}\right)}l^*(s).\end{equation}
To see this, we also observe that for $\Re
z<\si$
\[\frac{c_{i\mu}\GG\left(\frac{s-z}{2}\right)}{\GG\left(\frac{s-z}{2}+\frac{1}{2}-i\mu\right)}
=j^*(z-s),\] where
\[j(u)=(1-u^{-2})_+^{-1/2-i\mu}.\]
This is a special case of formula 3.251.1 from \cite{GR}.
\begin{Not}For $u\in\RR$, $u_+$ abbreviates $\max(0,u)$.\end{Not}
The required identity (\ref{eq98}) now reads
\[k^*(z)=\frac{1}{2\pi i}\int_{(\si)}j^*(z-s)l^*(s)\,ds,\quad 0<\Re z<\si<\si_0.\]
If we also assume that
\[\int_{(\si)}|l^*(s)|\,ds<\infty\]
(which will be the case, cf. (\ref{eq914})), then
a straightforward application of Fubini's theorem shows that the integral
evaluates the Mellin transform of $j(u)l(u)$ at $z$. Indeed,
\[\begin{split}(jl)^*(z)&=\int_0^\infty j(u)l(u)u^z\frac{\,du}{u} =\int_0^\infty
j(u)\left\{\frac{1}{2\pi
i}\int_{(\si)}l^*(s)u^{-s}\,ds\right\}u^z\frac{\,du}{u}\\\\
&=\frac{1}{2\pi i}\int_{(\si)}\left\{\int_0^\infty
j(u)u^{z-s}\frac{\,du}{u}\right\}l(s)\,ds=\frac{1}{2\pi
i}\int_{(\si)}j^*(z-s)l^*(s)\,ds.\end{split}\] As $k(u)$ and $j(u)l(u)$ are continuous,
(\ref{eq98}) is now equivalent to
\[k(u)=j(u)l(u).\]
If we use the fact and assumption that both $k(u)$ and $l(u)$ vanish for $u<1$,
this becomes 
\[l(u)=(1-u^{-2})_+^{1/2+i\mu}k(u).\] 
For $u>1$ the first factor can be expanded according to the binomial
theorem. The coefficients satisfy
\[\binom{\frac{1}{2}+i\mu}{j}\ll\frac{\GG\left(-\frac{1}{2}-i\mu+j\right)}{\GG(1+j)}
\ll(1+j)^{-3/2},\]
therefore Fubini's theorem yields
\begin{equation}\label{eq913}l^*(s)=
\sum_{j=0}^\infty\binom{\frac{1}{2}+i\mu}{j}(-1)^jk^*(s-2j).\end{equation}
From this representation and (\ref{eq220})--(\ref{eq912}) we can easily infer the bound
\begin{equation}\label{eq914}l^*(s)\ll_{\si,A}\bigl(1+|s|\bigr)^{-A},\quad0<\Re s<\si_0.\end{equation}

We found the recipe to construct a function $V^\jb(s)$ satisfying
(\ref{eq98}). First we determine $k^*$ according to
(\ref{eq912}). Then we define $l^*$ by (\ref{eq913}).
Finally, $V^\jb(s)$ is given by (\ref{eq911}). By our assumptions on $V^*(z)$, 
it is clear that $V^\jb(s)$ is analytic in
$0<\Re s<\si_0$, and 
from (\ref{eq914}) it also follows that in this region it satisfies 
a uniform upper bound
\[V^\jb(s)\ll_{\si,A}\bigl(1+|s|\bigr)^{-A-1/2}.\] Finally, (\ref{eq98}) implies (\ref{eq95})
for all points $u>0$, $u\neq 1$, because $V$ is continuous at all these
points by assumption. This completes the proof of the theorem.

\section{Bounds for hypergeometric functions}

In this section we prove Lemma~\ref{F}. 
For $0\leq u<\frac{1}{2}$ we use the representation
\begin{multline*}F\left(\tfrac{s}{2}+i\mu,\tfrac{1}{2}+i\mu;\tfrac{s}{2}+\tfrac{1}{2};u\right)=\\
\frac{\GG\left(\frac{s}{2}+\frac{1}{2}\right)}
{\GG\left(\frac{1}{2}+i\mu\right)\GG\left(\frac{s}{2}-i\mu\right)}
\int_0^1t^{-\frac{1}{2}+i\mu}(1-t)^{\frac{s}{2}-i\mu-1}(1-tu)^{-\frac{s}{2}-i\mu}\,dt.\end{multline*}
This identity is a special case of formula 9.111 from \cite{GR}.
It follows that
\[F\left(\tfrac{s}{2}+i\mu,\tfrac{1}{2}+i\mu;\tfrac{s}{2}+\tfrac{1}{2};u\right)
\ll_\si\left|\frac{\GG\left(\frac{s}{2}+\frac{1}{2}\right)}
{\GG\left(\frac{s}{2}-i\mu\right)}\right|
\int_0^1t^{-\frac{1}{2}}(1-t)^{-\frac{1}{2}}\,dt.\]
The integral on the right hand side is bounded, hence we have
\[F\left(\tfrac{s}{2}+i\mu,\tfrac{1}{2}+i\mu;\tfrac{s}{2}+\tfrac{1}{2};u\right)
\ll_{\si}|s|^{1/2}.\]

For the rest of this section we shall assume that 
$\frac{1}{2}\leq u<1$. We apply formula 9.131.1 from \cite{GR}:
\begin{equation}\label{eq922}
F\left(\tfrac{s}{2}+i\mu,\tfrac{1}{2}+i\mu;\tfrac{s}{2}+\tfrac{1}{2};u\right)
=(1-u)^{-\frac{s}{2}-i\mu}
F\left(\tfrac{s}{2}+i\mu,\tfrac{s}{2}-i\mu;\tfrac{s}{2}+\tfrac{1}{2};\tfrac{u}{u-1}\right).
\end{equation}
Note that here $\frac{u}{u-1}\leq-1$. We can express 
the hypergeometric function on the right hand side as a contour integral by
formula 9.113 from \cite{GR}:
\begin{multline*}
F\left(\tfrac{s}{2}+i\mu,\tfrac{s}{2}-i\mu;\tfrac{s}{2}+\tfrac{1}{2};\tfrac{u}{u-1}\right)=\\
\frac{1}{2\pi i}\int_{(-\ep)}
\frac{\GG\left(\frac{s}{2}+i\mu+w\right)
\GG\left(\frac{s}{2}-i\mu+w\right)
\GG\left(\frac{s}{2}+\frac{1}{2}\right)}
{\GG\left(\frac{s}{2}+i\mu\right)
\GG\left(\frac{s}{2}-i\mu\right)
\GG\left(\frac{s}{2}+\frac{1}{2}+w\right)}
\GG(-w)\left(\frac{u}{1-u}\right)^w\,dw.\end{multline*}
This formula is valid whenever
\[0<\ep<\frac{\si}{2}.\]

In order to estimate the integral efficiently, we shift the contour to
the line $\Re s=-\frac{\si}{2}-\ep$. This shift picks up the poles at
\[w=-\frac{s}{2}\pm i\mu.\]
To be precise, these are two simple poles when $\mu\neq 0$, and a double pole when $\mu=0$. In both cases we can write the result as
\begin{multline*}
F\left(\tfrac{s}{2}+i\mu,\tfrac{s}{2}-i\mu;\tfrac{s}{2}+\tfrac{1}{2};\tfrac{u}{u-1}\right)=
\sum_\pm d_{\pm i\mu}\frac{\GG\left(\frac{s}{2}+\frac{1}{2}\right)}{\GG\left(\frac{s}{2}\mp i\mu\right)}
\left(\frac{u}{1-u}\right)^{-\frac{s}{2}\mp i\mu}\\
+\frac{1}{2\pi i}\int_{\left(-\frac{\si}{2}-\ep\right)}
\frac{\GG\left(\frac{s}{2}+i\mu+w\right)
\GG\left(\frac{s}{2}-i\mu+w\right)
\GG\left(\frac{s}{2}+\frac{1}{2}\right)}
{\GG\left(\frac{s}{2}+i\mu\right)
\GG\left(\frac{s}{2}-i\mu\right)
\GG\left(\frac{s}{2}+\frac{1}{2}+w\right)}
\GG(-w)\left(\frac{u}{1-u}\right)^w\,dw,\end{multline*}
where $d_{i\mu}$ and $d_{-i\mu}$ are suitable constants. 
It follows from (\ref{eq922}) that
\begin{multline}\label{eq924}F\left(\tfrac{s}{2}+i\mu,\tfrac{1}{2}+i\mu;\tfrac{s}{2}+\tfrac{1}{2};u\right)
\ll_\si|s|^{1/2}\\
+\int_{\left(-\frac{\si}{2}-\ep\right)}
\left|\frac{\GG\left(\frac{s}{2}+i\mu+w\right)
\GG\left(\frac{s}{2}-i\mu+w\right)
\GG\left(\frac{s}{2}+\frac{1}{2}\right)}
{\GG\left(\frac{s}{2}+i\mu\right)
\GG\left(\frac{s}{2}-i\mu\right)
\GG\left(\frac{s}{2}+\frac{1}{2}+w\right)}
\GG(-w)\,dw\right|.\end{multline}

It remains to estimate the last integral. In the light of the uniform estimate
\[\left|\frac{\GG\left(\frac{s}{2}+i\mu+w\right)
\GG\left(\frac{s}{2}-i\mu+w\right)}
{\GG\left(\frac{s}{2}+i\mu\right)
\GG\left(\frac{s}{2}-i\mu\right)}\right|
\ll_{\si,\ep}\left|\frac{\GG^2\left(\frac{s}{2}+w\right)}
{\GG^2\left(\frac{s}{2}\right)}\right|,\quad\Re w=-\tfrac{\si}{2}-\ep,\]
we are left with estimating
\begin{equation}\label{eq923}\II=\int_{\left(-\frac{\si}{2}-\ep\right)}
\left|\frac{\GG^2\left(\frac{s}{2}+w\right)\GG\left(\frac{s}{2}+\frac{1}{2}\right)}
{\GG^2\left(\frac{s}{2}\right)\GG\left(\frac{s}{2}+\frac{1}{2}+w\right)}\GG(-w)\,dw\right|.\end{equation}
The value of the integral does not change when $s$ is replaced by $\bar s$, therefore we can
assume that $\Im s>0$. We split the integral into three parts.

\begin{Part} $\Im w>0$. In this segment Stirling's formula implies
\[\begin{split}
\frac{\GG^2\left(\frac{s}{2}+w\right)\GG\left(\frac{s}{2}+\frac{1}{2}\right)}
{\GG^2\left(\frac{s}{2}\right)\GG\left(\frac{s}{2}+\frac{1}{2}+w\right)}\GG(-w)
&\ll_{\si,\ep}e^{-\pi\Im w}|w|^{\frac{\si}{2}-\frac{1}{2}+\ep}
|s|^{1-\frac{\si}{2}}\left|\tfrac{s}{2}+w\right|^{-1-\ep}\\
&\ll_{\si,\ep}|s|^{1-\frac{\si}{2}}\left|\tfrac{s}{2}+w\right|^{-1-\ep}.
\end{split}\]
It follows that the total contribution to the integral (\ref{eq923}) is
\[\II_1\ll_{\si,\ep}|s|^{1-\frac{\si}{2}}|s|^{-\ep}\ll_{\si,\ep}|s|^{\frac{1}{2}}.\]
\end{Part}

\begin{Part} $0\geq\Im w\geq-\Im\tfrac{s}{2}$. In this segment Stirling's formula implies
\[\begin{split}
\frac{\GG^2\left(\frac{s}{2}+w\right)\GG\left(\frac{s}{2}+\frac{1}{2}\right)}
{\GG^2\left(\frac{s}{2}\right)\GG\left(\frac{s}{2}+\frac{1}{2}+w\right)}\GG(-w)
&\ll_{\si,\ep}|w|^{\frac{\si}{2}-\frac{1}{2}+\ep}
|s|^{1-\frac{\si}{2}}\left|\tfrac{s}{2}+w\right|^{-1-\ep}\\
&\ll_{\si,\ep}|s|^{\frac{1}{2}+\ep}\left|\tfrac{s}{2}+w\right|^{-1-\ep}.
\end{split}\]
It follows that the total contribution to the integral (\ref{eq923}) is
\[\II_2\ll_{\si,\ep}|s|^{\frac{1}{2}+\ep}.\]
\end{Part}

\begin{Part}$-\Im\tfrac{s}{2}>\Im w$. In this segment Stirling's formula implies
\[\begin{split}
\frac{\GG^2\left(\frac{s}{2}+w\right)\GG\left(\frac{s}{2}+\frac{1}{2}\right)}
{\GG^2\left(\frac{s}{2}\right)\GG\left(\frac{s}{2}+\frac{1}{2}+w\right)}\GG(-w)
&\ll_{\si,\ep}e^{\pi\Im\left(\frac{s}{2}+w\right)}|w|^{\frac{\si}{2}-\frac{1}{2}+\ep}
|s|^{1-\frac{\si}{2}}\left|\tfrac{s}{2}+w\right|^{-1-\ep}\\
&\ll_{\si,\ep}|s|^{\frac{1}{2}+\ep}
e^{\pi\Im\left(\frac{s}{2}+w\right)}\left|\tfrac{s}{2}+w\right|^{\frac{\si}{2}-\frac{3}{2}}.
\end{split}\]
It follows that the total contribution to the integral (\ref{eq923}) is
\[\II_3\ll_{\si,\ep}|s|^{\frac{1}{2}+\ep}.\]
\end{Part}

Altogether we can see that
\[\II=\II_1+\II_2+\II_3\ll_{\si,\ep}|s|^{1/2+\ep},\]
therefore (\ref{eq924}) and (\ref{eq923}) imply the required bound (\ref{eq97}).

The proof of Lemma~\ref{F} is complete.\qed

\end{document}